\documentclass[11pt]{amsart}
\usepackage{fullpage}
\allowdisplaybreaks
\usepackage{import}
\usepackage{hyperref}
\hypersetup{
    colorlinks=true,
    linkcolor=blue,
    filecolor=magenta,      
    urlcolor=cyan,
    citecolor=red,
    pdfpagemode=FullScreen,
    }
\usepackage{amsmath,amsfonts,amssymb}
\usepackage{mathrsfs}
\usepackage{verbatim}
\usepackage{multirow}

\usepackage{todonotes}
\newtheorem{theorem}{Theorem}[section]
 \newtheorem{corollary}[theorem]{Corollary}
 \newtheorem{lemma}[theorem]{Lemma}
 \newtheorem{proposition}[theorem]{Proposition}
 \theoremstyle{definition}
 \newtheorem{definition}[theorem]{Definition}
 
 \theoremstyle{remark}
 \newtheorem{remark}[theorem]{Remark}
  \newtheorem{ex}[theorem]{Example}
 \numberwithin{equation}{section}

\def \bC {\mathbb C}

\def \bN {\mathbb N}

\def \bR {\mathbb R}

\def \cC {\mathcal C}
\def \cD {\mathcal D}

\def \cR {\mathcal R}

\def \cW {\mathcal W}

\def \fg {\mathfrak g}

\def \id {\text{\rm I}}

\def \eps {\varepsilon}
\def \vol {{\rm vol}}
\def \pr  {{\rm pr}}
\def \Hom  {{\rm Hom}}
\def \IM  {{\rm Im}\, }

\def \Mat {{\rm Mat}} 
\def \spec  {{\rm sp}\, }
\def \diag  {{\rm diag} }

\begin{document}
\author[V. Fischer]{V\'eronique Fischer and Francesca Tripaldi}
\address[V. Fischer]%
{University of Bath, Department of Mathematical Sciences, Bath, BA2 7AY, UK.} 
\email{v.c.m.fischer@bath.ac.uk}

\address[F. Tripaldi]%
{Mathematisches Institut,
University of Bern,
Sidlerstrasse 5 3012 Bern, Switzerland.} 
\email{francesca.tripaldi@unibe.ch}

\title{An alternative construction of the  Rumin complex\\
on homogeneous nilpotent Lie groups}

\medskip



\maketitle

\begin{abstract}
In this paper, we consider the Rumin complex on homogenenous nilpotent Lie groups.  
We present an  alternative construction to the classical one on Carnot groups using  ideas from parabolic geometry. We also give the explicit computations for the Engel group with this 
approach.
\end{abstract}








\makeatletter
\renewcommand\l@subsection{\@tocline{2}{0pt}{3pc}{5pc}{}}
\makeatother

\tableofcontents

\section{Introduction}

What is now known as the Rumin complex was constructed by M. Rumin on contact manifolds in his PhD thesis \cite{Rumin1990,Rumin1994} around 1990.
The same name is also given to the complex he constructed on Carnot groups some years later \cite{Rumin1999,Rumin2000}, as the two objects coincide in the case of Heisenberg groups.
In this paper, we present an alternative approach to the Rumin complex on arbitrary Carnot groups.

\smallskip

Our setting is in fact slightly more general than Carnot groups. Using the vocabulary of G. Folland and E. Stein \cite{FollandStein}, we are considering homogeneous groups: these are the connected, simply connected nilpotent Lie groups with dilations (see Section \ref{subsec_homLgr}).
In our view, this is 
the most general setting where  the concept of \textit{homogeneous weights} coming from the group dilations can be applied  to construct the Rumin complex - either the original construction on Carnot groups, or the one presented here on homogeneous Lie groups.

The ideas behind the alternative construction we are presenting here come from parabolic geometry,
and in particular the Bernstein-{G}elfand-{G}elfand sequences, see \cite{CapSlovakSouvek,CalderbankDiemer}
and the 2017 preprint \cite{DaveHaller}.
They provide a robust and elegant approach that lends itself to future generalisations in more involved geometric contexts. We believe that they will allow for the construction of similar subcomplexes on 
equiregular sub-Riemannian manifolds
and more generally to filtered manifolds 
\cite{DaveHaller,tripaldi2020rumin}. 
Eventually, the treatment of singularities in the sub-Riemannian setting might also be possible after further progress.

M. Rumin's original motivations
 came from investigating the possible definitions of  curvature on contact manifolds starting from a  partial connection \cite{pansu1993differential}. His construction quickly found applications in analysis on complex hyperbolic spaces \cite{Rumin1994,RuminGAFA2000, vcap2020poisson}, as well as spectral geometry with the construction of hypoelliptic Laplace-type operators on contact manifolds
 \cite{Rumin+Seshadri2012, kitaoka2020analytic, pansu2019averages,BFP2,BFP3,BFP4}. 
 It should be noted that, starting from a different initial question, the Rumin complex had already been developed independently. Indeed,
the study of differential equations generalising Monge-Amp\`ere equations on contact manifolds led to the construction of an equivalent subcomplex by  
 Lychagin \cite{lychagin1979contact} in 1979 (and later expanded in \cite{lychagin1994non}), which was viewed as a variant of the
effective cohomology groups for contact manifolds \cite{bouche}.

In the past decade, the Rumin complex has become an essential tool in geometric analysis on contact manifolds and  on Carnot groups,
for instance regarding the study of quasi-isometry invariants \cite{pansu2017cup,Pansu+Rumin2018}, Lipschitz graphs \cite{FranchiSerapioni}, form-valued Partial Differential Equations \cite{FranchiTesi,baldi2015gagliardo, BFP1,FranchiMontefalconeSerra,baldi2022sobolev}, and more recently in developing a theory of sub-Riemannian currents on contact manifolds \cite{canarecci2021sub,vittone2022lipschitz}.
It is therefore of great interest to develop a theory behind the Rumin complex on Carnot groups that can enable us to construct an analogous subcomplex for arbitrary sub-Riemannian manifolds: this is the purpose of this paper.

\smallskip

The paper is organised as follows.
In Section \ref{sec_licohomG}, we recall the main definitions and properties of the left-invariant cohomology of an arbitrary Lie group $G$.
It is only starting from Section \ref{section 3} that we assume that $G$ is in addition nilpotent and homogeneous; in this context we recall Rumin's construction of the subcomplex that bears his name.
In Section \ref{sec_altconst}, we present the main result of the paper: an alternative construction of the Rumin complex on a nilpotent  homogeneous Lie group $G$.
Finally, we carry out the explicit computations behind this alternative construction in the case of the Engel group, i.e. the four-dimensional filiform nilpotent Lie group.

\subsection*{Acknowledgements}

V. Fischer acknowledges  the support of The Leverhulme Trust via Research Project Grant  2020-037, {\it Quantum limits for sub-elliptic operators}.

F. Tripaldi is supported by the Swiss National Science Foundation Grant nr 200020-191978, \textit{Analytic and geometric structures in singular spaces}.

\section{The 
left-invariant cohomology of a Lie group $G$}
\label{sec_licohomG}

In this section, 
we recall some well-known facts about the 
left-invariant cohomology of a given Lie group $G$, especially with regards to the Chevalley-Eilenberg differential $d_0$ and the Lie algebra cohomology of $G$. 
This is mainly  an opportunity to set some notation and conventions for the rest of the paper.

\subsection{The Chevalley-Eilenberg differential and co-differential}

\subsubsection{Generalities and notation}
Let us start with some generalities on the space
$\Omega^k M$ of smooth $k$-forms on a smooth manifold $M$.
We will write $\Omega^\bullet M= \oplus_k \Omega^k M$ and denote by 
$d^{(k)}:\Omega^k M \to \Omega^{k+1} M $ the action of the exterior derivative $d$ on forms of degree $k$.
For $k=0$, $\Omega^0 M $ is the space $C^\infty(M)$ of smooth functions on $M$, while for $k>0$ the space
 $\Omega^{k} M $ is a module over $\Omega^0 M =C^\infty(M)$. 
  For $k>0$,
 we have the following
explicit formula for the action of the exterior derivative $d$ on an arbitrary smooth $k$-form $\omega\in\Omega^kM$:
\begin{align}
\label{eq_def_dexplicit}
d^{(k)}\omega(V_0,\ldots,V_k)
&=
\sum_{i=0}^k (-1)^i V_i\left(\omega(V_0,\ldots,\hat V_i,\ldots, V_k) \right)
\\&\qquad +
\sum_{0\leq i<j\leq k}
(-1)^{i+j} \omega([V_i,V_j],V_0,\ldots,\hat V_i,\ldots,\hat V_j,\ldots, V_k)\,,\nonumber
\end{align}
where $V_0,\ldots,V_k\in\Gamma(TM)$ are vector fields on $M$, 
the bracket denotes the commutator bracket, and the hat an omission. 

\smallskip

Throughout the rest of the paper, 
we consider the case where the smooth manifold $M$ is a non-trivial Lie group $(G,\ast)$.
Note that nothing else will be assumed on $G$ throughout the present section. 

\smallskip

It is well-known that the exterior derivative $d$ commutes with the pullback of an arbitrary smooth function $f\colon G\to G$. In particular, this is true for left-translations of the group $L_g\colon G\to G$ with $ L_g(x)=g\ast x$, for a fixed $g\in G$. Therefore, we have that for any $g\in G$ and any smooth form $\omega\in\Omega^\bullet G$, $L_g^\sharp d\omega=d L_g^\sharp \omega$, and so $d$ may be viewed as a vector-valued left-invariant differential operator of degree 1. 

\subsubsection{Left-invariant forms}
	A $k$-form $\omega\in\Omega^kG$ is left-invariant if it is invariant under the pullback induced by the left-translations of the group, i.e. $L_g^\sharp \omega=\omega$ for any $g\in G$.
We denote the subspace of left-invariant $k$-forms by 
$$
\Omega_L^k G =\{\omega \in \Omega^k G  \mid L_g^\sharp \omega=\omega \, ,\ \forall g\in G\ \}  \ \subset \Omega^k G .
$$ 
 It is easy to check that $\Omega_L^0G =\bR $ and $d\big( \Omega^\bullet_LG\big)\subset\Omega^{\bullet}_LG$, since the pullback is compatible with the exterior derivative.

Recall that the subspace $\Gamma_L(TG)\subset\Gamma(TG)$ of left-invariant vector fields of $G$ is naturally isomorphic to its Lie algebra $\fg$, and so by duality the subspace of left-invariant 1-forms is naturally identified with 
$ \fg^*=\Hom(\fg,\bR)$. 
More generally, $\Omega_L^k G $ is naturally identified with $\Hom(\bigwedge ^k \fg ,\bR)$, for $ k=0,\ldots,\dim\mathfrak g$.

\subsubsection{The Chevalley-Eilenberg differential  $d_0$ }\label{def differential d_0}
As $d^{(k)}$ maps $\Omega_L^k G $ to $\Omega_L^{k+1} G $, 
and we identify each $\Omega_L^k G $ with $\Hom(\bigwedge^k \fg ,\bR)$, we obtain a map 
$d^{(k)}_0:\Hom(\bigwedge^k \fg ,\bR)\to \Hom(\bigwedge^{k+1} \fg ,\bR).$
	Such a map $d_0$ is obtained by simply restricting $d$ to $\Omega_L^\bullet  G \cong  \Hom(\bigwedge^\bullet \fg ,\bR)$, and it is commonly known as the Chevalley-Eilenberg differential.
	
The following properties 
for Chevalley-Eilenberg differential $d_0$ readily follow from the properties of the exterior differential $d$:
	\begin{itemize}

  \item[1.]  $d_0$ decomposes into linear maps
 $$
 d_0 =\sum_{k=-\infty}^\infty d_0^{(k)}
 \quad\mbox{with}\ d_0^{(k)}:\Hom(\bigwedge^k \fg ,\bR)\to \Hom(\bigwedge^{k+1} \fg ,\bR) , 
 $$
 with $d_0^{(k)} =0$ for $k\leq 0$ and $k\geq \dim \mathfrak g$;
 \item[2.] for $0<k< \dim \mathfrak g$, the explicit formula for each
	$d_0^{(k)}$ is given by  
\begin{equation}
\label{eq_d0}
d_0^{(k)}\omega\, (V_0,\ldots,V_k)
=
\sum_{0\leq i<j\leq k}
(-1)^{i+j} \omega([V_i,V_j],V_0,\ldots,\hat V_i,\ldots,\hat V_j,\ldots, V_k)
\end{equation}
for any $\omega\in \Hom(\bigwedge^k \fg ,\bR)$, and $V_0,\ldots, V_k\in \fg$;
 	\item[3.]  $(\Omega_L^\bullet G,d_0)$ is a co-chain complex, i.e. $d_0^{2}=0$ on $\Omega_L^\bullet G\cong \Hom(\bigwedge^\bullet\mathfrak{g},\mathbb R)$;
\item[4.] $d_0$ obeys the Leibniz rule, i.e. for $\alpha \in \Hom(\bigwedge^k \fg ,\bR)$ and $\beta\in \Hom(\bigwedge^\bullet \fg ,\bR)$, we have
$$
d_0(\alpha\wedge \beta) = (d_0\alpha)\wedge \beta + (-1)^k\alpha \wedge (d_0\beta).
$$ 

	\end{itemize}
Let us stress that the formula of $d_0$ presented in \eqref{eq_d0} is a direct consequence of applying the explicit formulation of the exterior derivative in
\eqref{eq_def_dexplicit} to left-invariant objects.
 
 \subsubsection{Scalar product on $\Hom(\bigwedge^\bullet \fg,\bR)$}
 \label{subsubsec_scalarproduct}
 
 We fix a basis $\lbrace X_1,\ldots, X_n\rbrace$ of $\fg$ (here we denote $n=\dim \fg$). 
 This yields a natural scalar product on $\Hom(\fg,\bR)$, which extends to  $\Hom(\bigwedge^k \fg,\bR)$  by imposing
$$
\langle \alpha_1\wedge \ldots\wedge \alpha_k,  \beta_1\wedge\ldots\wedge \beta_k\rangle_{\Hom(\bigwedge^k \fg,\bR)} := \det\left( \langle \alpha_i,\beta_j\rangle_{\Hom(\fg,\bR)}\right)_{1\leq i,j\leq k},
$$
and eventually to
$$
\Hom(\bigwedge^\bullet \fg,\bR) =\oplus_k^\perp\Hom(\bigwedge^k \fg,\bR)\ , 
\ \mbox{with}\ 
 \Hom(\bigwedge^k \fg,\bR) = 0 \ \mbox{for}\ k>n \ \mbox{or}\ k<0.
 $$
We may drop the indices of the scalar products $\Hom(\bigwedge^\bullet \fg,\bR)$ when the context is clear. 
Given a subspace $S$ of $  \Hom(\bigwedge^\bullet \fg,\bR)$, we will denote
the orthogonal projection onto $S$ by
$\pr_S$.

Given $\lbrace\theta^1,\ldots,\theta^n\rbrace\subset\Hom(\mathfrak g,\mathbb R)$ the dual basis of $\lbrace X_1,\ldots,X_n\rbrace$, i.e. $\theta^i(X_j)=\delta_{ij}$, then we set
$$
\theta^\emptyset:=1 \qquad\mbox{and}\qquad 
\vol:=\theta^{(1,\ldots,n)}:=\theta^1\wedge\ldots \wedge \theta^n,
$$
so that  $\Hom(\bigwedge^0 \fg,\bR) = \bR 1= \bR \theta^\emptyset$
and  $\Hom(\bigwedge^n \fg,\bR) = \bR \vol$.

It easy to check that
$\theta^{(i,j)}:=\theta^i\wedge \theta^j$, $1\leq i<j\leq n$
form a basis of $\bigwedge^2\mathfrak{g}^\ast=\Hom(\bigwedge^2 \fg,\bR)$, and
more generally,
	 the $k$-forms
\begin{equation}
\label{eq_thetaI}
\theta^I = \theta^{i_1}\wedge\ldots \wedge \theta^{i_k}\;,
\mbox{ with increasing indices} \ I=(i_1,\ldots,i_k), \ 1\leq i_1<\cdots <i_k\leq n, 
\end{equation}
form an orthonormal basis of $\bigwedge^k\mathfrak{g}^\ast=\Hom(\bigwedge^k \fg,\bR)$ for any $k=1,\ldots, n$.

\subsubsection{The Chevalley-Eilenberg co-differential  $d_0^t $ }\label{def adjoint of d_0}
The Chevalley-Eilenberg co-differential
is the transpose of $d_0$ for the scalar product on $\Hom(\bigwedge^\bullet \fg,\bR)$ defined above, and we will denote it by $d_0^t$.
In other words,  $d_0^t$ is the linear map acting on $\Hom(\bigwedge^\bullet \fg,\bR)$ such that 
$$
\langle d_0 \alpha , \beta\rangle _{\Hom(\bigwedge^\bullet \fg,\bR) }
=\langle \alpha , d_0^t \beta\rangle_{\Hom(\bigwedge^\bullet \fg,\bR) }\;, 
\quad\forall\ \alpha,\beta\in \Hom(\bigwedge^\bullet \fg,\bR). 
$$
The properties of $d_0$ imply that its transpose decomposes as 
$$
d_0^t =  \sum_{k=-\infty}^\infty d_0^{(k,t)}\ , 
\ 
d_0^{(k,t)}:\Hom(\bigwedge^{k+1} \fg,\bR) \to \Hom(\bigwedge^{k} \fg,\bR)\ ,
$$ 
and  $d_0^{(k,t)}=0$ for $k\leq 0$ and $k\ge n$. 
Moreover,   $(d_0^{t})^2=0$.

 \subsubsection{The Hodge star operator on $\Hom(\bigwedge^\bullet \fg,\bR)$}
 \label{subsubsec_star0}

The Hodge star operator is the map acting on $\Hom(\bigwedge^\bullet \fg,\bR)$ via 
$$
\star : \Hom(\bigwedge^k \fg,\bR) \overset {\cong} \longrightarrow \Hom(\bigwedge^{n-k} \fg,\bR), 
\quad 
\alpha\wedge \star\beta := \langle\alpha,\beta\rangle \vol\quad\forall \,\alpha,\beta \in \Hom(\bigwedge^k \fg,\bR) .
$$
From the definition it easily follows that
$$
\star 1=\vol, \qquad \star \vol =1, \qquad 
\alpha \wedge \star \beta = \beta \wedge \star \alpha,
$$
and
$$
\star \theta^i = (-1)^{i+1}\theta^1\wedge\ldots\wedge \widehat \theta^i\wedge \ldots\wedge \theta^n\ .
$$
More generally, the Hodge star operator on $\theta^I$ with increasing indices $I=(i_1,\ldots,i_k)$
is given by
$$
\star (\theta^I)
=
(-1)^{\sigma(I)}
\theta^{\bar I}
$$
where $\bar I = (\bar i_1,\ldots,\bar i_{n-k})$ with $\bar i_1<\bar i_2<\ldots <\bar i_{n-k}$ complements $I$, and $(-1)^{\sigma(I)}$ is the sign of the permutation $\sigma(I)=i_1\cdots i_k\,\bar i_1\cdots\bar i_{n-k}$.

This characterises the Hodge star operator on a given basis of $ \Hom(\bigwedge^\bullet \fg,\bR)$, and
hence by linearity on the whole $ \Hom(\bigwedge^\bullet \fg,\bR)$.

Moreover, the previous properties imply readily that  for any $\alpha,\beta\in \Hom(\bigwedge^k \fg,\bR)$, we have
$$
\star\star \beta = (-1)^{k(n-k)}\beta
\qquad\mbox{and}\qquad\langle\star \alpha,\star \beta\rangle= \langle\alpha,\beta\rangle\,,
$$
and that  the transpose $d_0^t$ of $d_0$ satisfies
$$
d_0^{(k,t)}  = (-1)^{kn +1} \star d_0^{(n-k-1)} \star\,.
$$

\subsection{The subspaces $E_0$ and $F_0$}
\label{subsecE0F0}

For each $k=0,\ldots,n$, let us consider the subspace of $\Hom(\bigwedge^k \fg,\bR)$ defined by
$$
F_0^{k} := \IM d_0^{(k-1)} + \IM d_0^{(k,t)}\subset\Hom(\bigwedge^k \fg,\bR)\,,
$$
where the maps
\begin{align*}
    d_0^{(k-1)}:\Hom(\bigwedge^{k-1} \fg ,\bR)\to \Hom(\bigwedge^{k} \fg ,\bR)\ \text{ and }\ d_0^{(k,t)}:\Hom(\bigwedge^{k+1} \fg,\bR) \to \Hom(\bigwedge^{k} \fg,\bR)
\end{align*}
were the maps defined in Subsection \ref{def differential d_0} and \ref{def adjoint of d_0}.

\begin{lemma}
\label{lem_Hd0F0} 
The complex $(F_0^\bullet,d_0)$ is  acyclic, i.e. $d_0(F_0^\bullet)\subset F_0^\bullet$ and the resulting cohomology is trivial.
\end{lemma}
\begin{proof}
Given that $(\Hom(\bigwedge^\bullet\mathfrak g,\mathbb R),d_0)$ is a complex, i.e. $d_0^2=0$, and applying the definition of $d_0^{t}$, we obtain
\begin{align*}
  d_0^{(k)} (F_0^{k}) & = d_0^{(k)} (\IM d_0^{(k,t)}) = d_0^{(k)} (\ker d_0^{(k)})^\perp = \IM d_0^{(k)}\  \subset F_0^{k+1} \,.
\end{align*}
Moreover
\begin{align*}
    F_0^k \cap \ker d_0^{(k)} =&
\IM d_0^{(k-1)} \cap \ker d_0^{(k)} + \IM d_0^{(k,t)} \cap \ker d_0^{(k)} 
\\=&\IM d_0^{(k-1)} + (\ker d_0^{(k)})^\perp \cap \ker d_0^{(k)} 
=\IM d_0^{(k-1)}\,\text{ and}\\
F_0^{k}\cap \IM d_0^{(k-1)}=&\IM d_0^{(k-1)} + (\ker d_0^{(k)})^\perp \cap \IM d_0^{(k-1)} 
=\IM d_0^{(k-1)}\,,
\end{align*}
so that
\begin{align*}
    H^k(F_0,d_0)=\frac{\ker d_0^{(k)}\colon F_0^k\to F_0^{k+1}}{\IM d_0^{(k-1)}\colon F_0^{k-1}\to F_0^k}=\lbrace 0\rbrace\ ,\ k=0,\ldots,n\,.
\end{align*}
\end{proof}

\noindent\textit{Convention for cohomology spaces.} In the proof above, we have used the following notation that we will also apply  
throughout this paper when computing the cohomology of  complexes with various exterior differentials. The $k^{th}$-cohomology space of a given complex $(\cC^\bullet,d_\cC)$ is denoted by:
\begin{align*}
    H^k(\cC,d_\cC):=\frac{\ker d_\cC\colon \cC^k\to \cC^{k+1}}{\IM d_\cC\colon \cC^{k-1}\to \cC^{k}}\,.
\end{align*}

Given Lemma \ref{lem_Hd0F0}, 
it makes sense to focus our attention on the orthogonal complement of $F_0^k$ in 
$\Hom(\bigwedge^k \fg,\bR)$:
$$
E_0^k:=(F_0^{k})^\perp = 
(\IM d_0^{(k-1)})^\perp \cap (\IM d_0^{(k,t)})^\perp
= \ker d_0^{(k-1,t)} \cap \ker d_0^{(k)} \subset\Hom(\bigwedge^k\mathfrak g,\mathbb R)\,.
$$

In order to simplify the notation, for the remainder of the paper we will drop the index $k$, unless we want to highlight the degree of the forms the given operator is acting on.

The main result of this section is the following.

\begin{proposition}
\label{prop_prE0}
The orthogonal projection $\pr_{E_0}$ onto $E_0\subset\Hom(\bigwedge^\bullet\mathfrak g,\mathbb R)$ may be constructed in the two following equivalent ways:
\begin{description}
\item[Rumin's construction]
as above, the partial inverse $d_0^{-1}$ of $d_0$ is defined as
$d_0^{-1} := (d_0)^{-1} \pr_{\IM d_0}$; the operator
$r_0 := \id -d_0^{-1}d_0 -d_0 d_0^{-1}$
coincides with $\pr_{E_0}$ and satisfies 
\begin{equation}
\label{eq_r0d0}
 	d_0^{-1}r_0 =r_0d_0^{-1} =0\quad,\quad
 	d_0r_0=r_0 d_0=0\,.
\end{equation}
\item[Construction with $\Box_0$]
The operator $\Box_0 := d_0d_0^t + d_0^t d_0$
is symmetric on $\Hom(\bigwedge^\bullet \fg,\bR)$. Moreover, the spectral projector $\Pi_0$ onto the 0-eigenspace of $\Box_0$ coincides with $\pr_{E_0}$ and satisfies 
\begin{equation}
\label{eq_Pi0d0}
 	d_0^{t} \Pi_0 =\Pi_0d_0^{t} =0\quad,\quad
 	d_0\Pi_0 =\Pi_0 d_0=0\,.
\end{equation}
\end{description}

In both constructions, we obtain that this projection $\pr_{E_0}=\Pi_0=r_0$ commutes with the Hodge star operator and preserves the degree of forms:
$$
\pr_{E_0^{n-k}}\, \star = \star\, \pr_{E_0^k} 
\qquad\mbox{and}\qquad
\pr_{E_0^k}(\Hom (\bigwedge^k \fg,\bR)) \subset \Hom (\bigwedge^k \fg,\bR).
$$
\end{proposition}

The next two sections are devoted to the proof of Proposition \ref{prop_prE0}:
we first sketch out the well known construction by Rumin in Section \ref{subsec_RPi0}, and then give the alternative construction with $\Box_0$ in 
Section \ref{subsec_altPi0}.

\subsection{Rumin's construction of $\pr_{E_0}$}
\label{subsec_RPi0}

 \subsubsection{The partial inverse $d_0^{-1}$} 
\label{subsec_d0inv}

The map $d_0$ is a bijection from $(\ker d_0)^\perp = \IM d_0^{t}$ onto $\IM d_0$.
Following Rumin's construction, it is customary to use the shorthand $d_0^{-1}$ for the linear map given by $(d_0)^{-1} \pr_{\IM d_0}$, so that
\begin{align*}
    d_0^{-1} = \sum_{k=-\infty}^\infty d_0^{(k,-1)}\ ,\ \ d_0^{(k,-1)} := (d_0^{(k)})^{-1} \pr_{\IM d_0^{(k)}} : \Hom(\bigwedge^{k+1} \fg,\bR) \longrightarrow 
\Hom(\bigwedge^k \fg,\bR)
\end{align*}
with $d_0^{(k,-1)}=0$ when $k\leq 0$ or $k\geq n$, and for $0<k<n$ we have $\IM d_0^{(k,-1)}  \subset
(\ker d_0^{(k)})^\perp $.

By construction, one can check that
\begin{itemize}
    \item[i.] $\ker d_0^{-1}=\ker d_0^{t}=(\IM d_0)^\perp$;
    \item[ii.] $d_0^{-1}d_0=\pr_{\IM d_0^{t}}$;
    \item[iii.] $d_0d_0^{-1}=\pr_{\IM d_0}$;
    \item[iv.] $(d_0^{-1})^2=0$.
\end{itemize}

\subsubsection{Matrix representations of the maps $d_0$, $d_0^t$, and $d_0^{-1}$} It can be useful to have explicit matrix representations $\Mat_{\mathcal{D}}^{\mathcal{R}}$ of the linear maps $d_0$, $d_0^t$ and $d_0^{-1}$, once some particular orthonormal bases $\mathcal{D}$ and $\mathcal{R}$ for the domain and range spaces respectively have been chosen. 

In fact, if we take $e$ and $f$ to be orthonormal bases of $\ker d_0^{(k)}$ and $(\ker d_0^{(k)})^\perp = \IM d_0^{(k,t)}$ respectively, while  
$h$ and $g$ are  orthonormal bases of $\IM d_0^{(k)}$ and $(\IM d_0^{(k)})^\perp = \ker d_0^{(k,t)}$ respectively, then $\mathcal{D}=(e, f)$ is an orthonormal basis of $\Hom(\bigwedge^k\mathfrak{g},\bR)$ and $\mathcal{R}=(h,g)$ is an orthonormal basis of $\Hom(\bigwedge^{k+1}\mathfrak{g},\bR)$. Moreover,
for this choice of bases, we obtain a particularly insightful matrix representation of the linear maps $d_0^{(k)}$, $d_0^{(k,t)}$ and $d_0^{(k,-1)}$. Indeed,
$$
\Mat_{\cD}^{\cR}\big(d_0^{(k)}\big)
=\begin{pmatrix}
	0&0\\0&A
\end{pmatrix},
$$
where $A$ is an invertible square matrix, and 
$$
\Mat_{\cR}^{\cD}\big(d_0^{(k,-1)}\big)
=\begin{pmatrix}
	0&0\\0&A^{-1}
\end{pmatrix}, 
\qquad 
\Mat_{\cR}^{\cD}\big(d_0^{(k,t)}\big)
=\begin{pmatrix}
	0&0\\0&A^{^t}
\end{pmatrix}.
$$

\subsubsection{Properties of the operator $r_0$}

By construction, 
the operator 
 \begin{align*}
     r_0^{(k)}:=\id-d_0^{(k,-1)}d_0^{(k)} -d_0^{(k-1)}d_0^{(k-1,-1)}
 \end{align*}
 preserves the degrees of forms. Moreover, by applying the properties of $d_0$ and $d_0^{-1}$ and carrying out the computations, it is easy to check that
$r_0^2=r_0$ and \eqref{eq_r0d0} hold.
Hence, the map $r_0$ is a projection onto 
 $$
 \IM r_0=  
 \ker d_0 \cap \ker d_0^{-1} =\ker d_0 \cap \ker d_0^{t} =E_0\,,
 $$
 along 
$$
\ker r_0
=
\IM d_0^{-1} +\IM d_0
=
\IM d_0^{t} +\IM d_0 =F_0\,.
$$
In other words, $r_0=\pr_{E_0}$.

Finally, these descriptions together with the fact that $d_0^{(k,t)}  = (-1)^{kn+1} \star d_0^{(n-k-1)} \star$
 imply that the operator $r_0$ commutes with the Hodge star operator $\star$.
 This concludes the sketch of the proof of \eqref{eq_r0d0} behind Rumin's construction in Proposition \ref{prop_prE0}.

\subsection{The alternative construction of $\pr_{E_0}$ with $\Box_0$}
\label{subsec_altPi0}

Consider the box operator associated with $d_0$:
$$
\Box_0 := d_0 d_0^t +d_0^t d_0.
$$
This is an operator acting on $\Hom(\bigwedge^\bullet\fg ,\bR)$ that respects the degrees of forms:
$$
\Box_0 = \sum_{k=-\infty}^\infty \Box_0^{(k)}, \ \mbox{with} \ 
\Box_0^{(k)}:=d_0^{(k-1)} d_0^{(k-1,t)} +d_0^{(k,t)} d_0^{(k)}  : \Hom(\bigwedge^k \fg,\bR)\to \Hom(\bigwedge^{k} \fg,\bR).
$$
The expression of $\Box_0^{(k)}$ for certain choices of $k$ can be explicitly computed:
$$
\Box_0^{(0)}=d_0^{(0,t)}d_0^{(0)}=0\ , \ \Box_0^{(1)}= d_0^{(1,t)} d_0^{(1)}\ ,\ \Box_0^{(n)}=d_0^{(n-1)} d_0^{(n-1,t)} \ ,
$$
and $\Box_0^{(k)}=0$ for $k< 0$ and $k>n$.

This operator $\Box_0$ is symmetric, given that for any $\alpha,\beta\in\Hom(\bigwedge^\bullet\mathfrak g,\bR)$ we have
$$
\langle\Box_0\alpha,\beta\rangle=\langle \alpha,\Box_0\beta\rangle\,.
$$
Moreover, its kernel is
$$
\ker \Box_0 = \ker d_0 \cap \ker d_0^{t}\ , 
$$
since clearly $\ker \Box_0 \supset \ker d_0 \cap \ker d_0^{t}$, 
while conversely for any $\alpha\in\ker \Box_0$
$$
0=\langle\Box_0 \alpha,\alpha\rangle=\langle d_0 d_0^{t} +d_0^{t} d_0\alpha,\alpha\rangle =\vert d_0 \alpha\vert^2 + \vert d_0^{t} \alpha\vert^2, $$
and therefore $ \alpha\in \ker d_0 \cap \ker d_0^{t}$.

Note that this final computation also shows that $\Box_0$ is non-negative, but we will not need this remark in this paper. 

For the sake of clarity, we include the following lemma, which states the spectral theorem for finite-dimensional symmetric map and some corollaries with the Cauchy integral formula. 
\begin{lemma}
\label{lem_matrix}
	Let $M$ be a symmetric linear map on a Euclidean space (i.e. a finite dimensional real vector space equipped with a scalar product).
\begin{itemize}
	\item[1.] The map $M$ is diagonalisable in an orthonormal basis.
	In other words, 
	$M=\sum_\lambda \lambda \Pi_\lambda$, where  $\Pi_\lambda$ is the orthogonal projection onto the $\lambda$-eigenspace $\ker (M-\lambda)$.
	The sum is taken over the spectrum $\spec(M)$ of $M$, that is, the finite set of eigenvalues of $M$. 
	\item[2.] 	Each spectral projection $\Pi_\lambda$ may be recovered from the resolvant $(z-M)^{-1}$, for
	$z\in \bC \setminus
	\spec (M)$, using the Cauchy integral:
	$$
	\Pi_\lambda = \frac1{2\pi i}\oint_\gamma (z-M)^{-1} dz \,,
	$$
	where $\gamma$ can be taken as a circle around $\lambda$ with radius small enough (more generally, $\gamma$ is a closed curve
	lying in $\bC \setminus
	\spec (M)$ with winding number 1 about $\lambda$ and 0 around all the other eigenvalues).  
\end{itemize}
\end{lemma}

Since each $\Box_0$ is a symmetric operator on $\Hom(\bigwedge^\bullet\mathfrak g,\bR)$, point 1. Lemma \ref{lem_matrix} yields the following spectral decomposition 
$$
\Box_0 = \sum_{\lambda \in \spec(\Box_0)} \lambda \Pi_{\lambda}\ ,
$$ 
where $\Pi_{\lambda}$ denotes the orthogonal projection onto the eigenspace for each eigenvalue $\lambda$. 

For $\lambda=0$, $\Pi_0$  is the orthogonal projection onto $\ker \Box_0= \ker d_0 \cap \ker d_0^{t} =E_0$, 
and therefore it coincides with $\pr_{E_0}$.
In particular,  $d_0 \Pi_0  = d_0^{t} \Pi_0=0$.

Since $d_0^2=(d_0^{t})^2=0$, we readily have that both $d_0$ and $d_0^{t}$ commute with $\Box_0$ and therefore also with $\Pi_0$ by point 2. of Lemma \ref{lem_matrix}.
Hence \eqref{eq_Pi0d0} holds.

Since $\Box_0$ respects the degrees of forms, so do its orthogonal projections again by point 2. of  Lemma \ref{lem_matrix}:
$$\Pi_{\lambda} = \sum_{k=-\infty}^\infty \Pi_{\lambda}^{(k)}\ ,
\ \mbox{ with }\
\Pi_{\lambda}^{(k)} :  \Hom(\bigwedge^k \fg,\bR)\to \Hom(\bigwedge^{k} \fg,\bR)\ .
$$
Finally, the properties of the Hodge star product on $\Hom(\bigwedge^k\mathfrak g ,\bR)$ imply that
$$
\star\, \Box_0^{(k)} \star  = (-1)^{k(n-k)} \Box_0^{(n-k)},
 $$
and again by point 2. of Lemma \ref{lem_matrix},
we obtain  $\star\, \Pi_{\lambda^{(k)}}^{(k)}\star =  \Pi_{\mu^{(n-k)}}^{(n-k)} $ with $\lambda^{} = (-1)^{k(n-k)}\mu$, $\lambda$ and $\mu$ being eigenvalues of $\Box_0$. 
In particular,  
$\Pi_0$ commutes with $\star$. 
 This concludes the  proof of  Proposition \ref{prop_prE0}.

 \subsection{Natural extensions to $\Omega^\bullet G$}
\label{subsec_extOmega}
Since $\Omega^\bullet G$  is a $C^\infty(G)$-module generated by  $\Omega_L^\bullet G$, 
many objects extend naturally from $\Omega_L^\bullet G \cong \Hom(\bigwedge^\bullet \fg,\mathbb R)$ to  $\Omega^\bullet G$, and we will mostly keep the same notation for their extensions. 

This is the case for the natural extension of the Hodge star operator  
$$
\star :\Omega^kG \overset {\cong} \longrightarrow \Omega^{n-k}G\ .
$$
Also for the elements of $\Omega^\bullet G$ we will still use the same notation we used for the $k$-covectors $\Hom(\bigwedge^\bullet \fg,\bR)$, which can be viewed as left invariant forms in $\Omega^\bullet G$. For instance,  $\theta^1, \ldots,\theta^n$ are smooth 1-forms, $\vol = \theta^1 \wedge \ldots\wedge \theta^n$ is a smooth $n$-form, and more in general $\theta^I$ given in \eqref{eq_thetaI} will be a basis for $\Omega^kG$.

The scalar product on $\Hom(\bigwedge^\bullet \fg,\bR)\cong \Omega_L^\bullet G$ extends naturally to the space of smooth forms $\Omega^\bullet G$. 
In particular, when considering the subspace $\Omega^\bullet_cG$ of forms in $\Omega^\bullet G$, one can introduce a scalar product, the so-called $L^2$-inner product on forms, defined as 
$$
\langle\alpha,\beta\rangle_{L^2(\Omega^\bullet G)}
:
= \int_G \alpha \wedge \star \beta,
\quad 
\alpha, \beta\in \Omega^\bullet_cG \, . 
$$
Using this scalar product on $\Omega^\bullet_cG$, one can define the formal transpose $d^t$ of $d$, often called the co-differential. More precisely,  $d^{t}$
is the linear map on $\Omega^\bullet G$ characterised by 
$$
\langle d \alpha , \beta\rangle_{L^2(\Omega^\bullet G)} = \langle\alpha, d^{t} \beta\rangle_{L^2(\Omega^\bullet G)}\;\ ,
\quad 
 \ \alpha, \beta\in \Omega^\bullet_cG  \,. 
$$

Routine checks show that $d^t$ is a differential operator satisfying
$$
d^t=\sum_{k=0}^{n-1} d^{(k,t)}, \qquad 
 d^{(k,t)}:\Omega^{k+1} G \to \Omega^k G, 
 \quad\mbox{with}\quad
d^{(k,t)}  = (-1)^{kn+1} \star d^{(n-k-1)} \star. 
$$ 

The operators $d_0$, $d_0^{t}$, $d_0^{-1}$, $\Box_0$, $\pr_S$ for any subspace $S\subset \Hom(\bigwedge^\bullet \fg,\bR)$ and in particular $\pr_{E_0}=\Pi_0$, extend naturally into algebraic linear operators on $\Omega^\bullet G$, and we will keep the same notation for these extensions as well as for the subspaces
$$
E_0:=\IM \Pi_0=\ker d_0 \cap \ker d_0^{t} \subset \Omega^\bullet G
\ \mbox{ and }\ 
F_0:=\ker \Pi_0=\IM d_0 + \IM d_0^{t} \subset \Omega^\bullet G\ .
$$
Just like before, the algebraic operators $d_0$, $d_0^t$, and $d_0^{-1}$ satisfy
$$
d_0^2=0\ ,\ (d_0^{-1})^2=0\ ,\ (d_0^t)^2=0\ \text{ and }
d_0^{(k,t)}  = (-1)^{kn+1} \star d_0^{(n-k-1)} \star\ , 
$$
and 
$$
\langle d_0 \alpha , \beta\rangle_{L^2(\Omega^\bullet G)} = \langle\alpha, d_0^t \beta\rangle_{L^2(\Omega^\bullet G)}\;\ ,
\quad 
 \alpha, \beta\in \Omega^\bullet_cG \, . 
$$
The operators $\Box_0$ and $\Pi_0$ respect the degrees of forms. The operator $\Pi_0$ is  the projection 
onto $E_0\subset \Omega^\bullet G$ along 
$F_0\subset \Omega^\bullet G$, 
and satisfies $\Pi_0^2=\Pi_0^t =\Pi_0$ and $\star \Pi_0 = \Pi_0 \star$. 
Here, $\Pi_0^t$ is the linear map on $\Omega^\bullet G$ characterised by 
$$
\langle \Pi_0 \alpha , \beta\rangle_{L^2(\Omega^\bullet G)} = \langle\alpha, \Pi_0^t \beta\rangle_{L^2(\Omega^\bullet G)},
\quad 
\alpha, \beta\in \Omega^\bullet_cG \, . 
$$

\section{The  Rumin complex  on homogeneous Lie groups}\label{section 3}

In this section, 
we assume that the Lie group $G$ is homogeneous in the sense of G. Folland and E. Stein \cite{FollandStein}. 
After explaining this notion in more detail in Section \ref{subsec_homLgr},
we will recall Rumin's construction of the complex $(E_0^\bullet,d_c)$ in this setting and eventually give an alternative but equivalent construction.

\subsection{The homogeneous Lie group $G$ and its Lie algebra $\fg$}
\label{subsec_homLgr}

Let $G$ be a homogeneous Lie group, i.e.  $G$ is a  connected simply connected Lie group whose real finite dimensional Lie algebra
$\fg$ admits a family of dilations $\delta_r$, $r>0$. This means that the $\delta_r$'s are morphisms of the Lie algebra $\fg$ of the form 
$\delta_r = {\rm Exp} (A \ln r)$ where $A$ is a diagonalisable operator on $\fg$ with positive eigenvalues, and  ${\rm Exp}$ is the exponential on the space of morphisms of the vector space $\fg$.

\subsubsection{The weights in $\fg$}\label{def weights lie algebra}
Let $X_1,\ldots, X_n$ be a basis of eigenvectors of $A$. 
In this basis, $A$ is represented by the diagonal matrix $\Mat(A)= {  \rm diag} (\upsilon_1,\ldots,\upsilon_n)$ with positive eigenvalues $\upsilon_j>0$, $j=1,\ldots,n$.
Hence, $\delta_r$  takes the matrix form $\Mat( \delta_r)= { \rm diag} (r^{\upsilon_1},\ldots,r^{\upsilon_n})$.
The eigenvalues $\upsilon_1,\ldots,\upsilon_n$, are called the weights of the dilations, 
and we denote the set of weights by 
\begin{align}\label{weights upsilon of g}
   \cW(\fg) := \{\upsilon_1,\ldots,\upsilon_n\}\,. 
\end{align}
We may assume that the eigenvalues are increasingly ordered:
$$
0 <\upsilon_1\leq \ldots\leq \upsilon_n\,.
$$
Let $w_1,\ldots,w_s \in \bR$ denote the weights listed without multiplicity:
$$
\cW(\fg) = \{\upsilon_1,\ldots,\upsilon_n\}
=\{w_1, \ldots ,w_s\},\quad\mbox{with}\quad
0<w_1< \ldots <w_s\,.
$$
For each weight $w_j$, we denote by $\fg_{w_j}$ the $w_j$-eigenspace for $A$. 
The Lie algebra decomposes into a direct sum 
$$
\fg = \oplus_{j=1}^s \fg_{w_j}
\qquad\mbox{satisfying}\qquad
[\fg_{w_i},\fg_{w_j}] \subseteq \fg_{w_i +w_j}, \
1\leq i,j\leq s\,.
$$
As a direct consequence, it follows that  $\fg$ is  nilpotent. 

Moreover, the Lie group $G$ being nilpotent, connected, and simply connected, implies that its exponential mapping $\exp: \fg \to G$ is a bijection, and furthermore a global diffeomorphism. 
This allows us to extend the dilations to the group.
Keeping the same notation for the dilations on the group, 
they are the maps $\delta_r$, with $r>0$, on $G$ characterised by
$\delta_r \exp V = \exp \delta_r V$,  $V\in \fg$.

\subsubsection{Carnot, stratified and graded groups as homogeneous groups}
Let us recall some particular cases of homogeneous groups often considered in  other  publications on the Rumin complex, 
potentially with a different vocabulary from 
G.  Folland and E. Stein \cite{FollandStein}.

When the weights of the dilations may be assumed to be integers \cite[Section 3.1]{FischerRuz}, the group is said to be \emph{graded}.
If in addition  $\fg_{w_1}$ with $w_1=1$ generates the whole Lie algebra $\fg$,  the Lie group $G$ is called \emph{stratified};
it is also known as \emph{Carnot} when a scalar product has been fixed on $\fg_1$ (this can always be assumed).
The main example of Carnot groups, especially in relation to the Rumin complex, is the Heisenberg group, see e.g. 
\cite{BFP2,FranchiMontefalconeSerra,FranchiSerapioni,FranchiTesi}.

Schematically, we have the following implications
$$
\mbox{stratified}\ 
\Longrightarrow\ 
\mbox{graded}\ 
\Longrightarrow\ 
\mbox{homogeneous}\ 
\Longrightarrow\ 
\mbox{nilpotent}.
$$
However, none of them can be reversed \cite[Section 3.1]{FischerRuz}: not all (connected simply connected) nilpotent Lie groups can be equipped with a homogeneous structure,
there are homogeneous groups that do not admit any gradation, and finally it is not always possible to find a stratification on every graded group.

\subsection{Weights of forms}
The dilations $\delta_r$, $r>0$, extend naturally to 
$\Hom(\bigwedge^k \fg,\bR)\cong \Omega_L^k G $ and then to $\Omega^k G$ for any $k\in \bN_0$ via 
$$
(\delta_r\theta)(x;V_1,\ldots,V_k) =\theta(x;\delta_r V_1,\ldots,\delta_r V_k)\ ,\ \
 \theta\in \Omega^k G \ ,\ x\in G\ , \ V_1,\ldots,V_k\in \Gamma (TG)|_x\cong \fg\,.
$$
We keep the same notation $\delta_r$ for these extensions.
Routine checks show that they respect the wedge product:
\begin{equation}
\label{eq_deltarwedge}
     \forall \alpha_1,\alpha_2 \in \Omega^\bullet G
     \qquad
     \delta_r (\alpha_1 \wedge \alpha_2) = (\delta_r \alpha_1) \wedge (\delta_r \alpha_2).
\end{equation}

\begin{definition}
We say that a form $\theta\in \Omega^\bullet G $ is \emph{homogeneous} when there exists  $w\in \bR$ such that $\delta_r \theta = r^w \theta$ for all $r>0$.
The number $w$ is then called the \emph{weight} of $\theta.$
\end{definition}

\begin{ex}
\label{ex_theta0}
The 0-form $\theta^\emptyset =1$  is homogeneous with weight 0.
\end{ex}

\begin{ex}
\label{ex_thetaiweight}
For each $i=1,\ldots,n$,
 the left-invariant 1-form $\theta^i\in\Omega^1G$ given via $X^\ast_i\in\mathfrak g^\ast$ is homogeneous, and its  weight is  $\upsilon_i$, the weight of $X_i$ (see \eqref{weights upsilon of g}).
 Indeed, we obtain the formula 
 $$
 (\delta_r \theta_i)=r^{\upsilon_i} \theta_i, \qquad r>0, \ i=1,\ldots, n,
 $$
 from the following computation
 $$
 (\delta_r \theta_i) (X_j) = \theta_i(\delta_r X_j) = r^{\upsilon_j} \theta_i(X_j) 
 = \left\{\begin{array}{ll}
 r^{\upsilon_i} \ & \mbox{if} \ i=j\,,\\
 0 \ & \mbox{if} \ i\neq j\,.
 \end{array}\right.
 $$
\end{ex}

Since the dilations respect the wedge product (see \eqref{eq_deltarwedge}), we check readily the following property.
\begin{lemma}
\label{lem_wedgeweight}
    If the forms $\alpha_1, \alpha_2\in \Omega^\bullet G$ are homogeneous with weights $w_1$ and $w_2$ respectively, then $\alpha_1\wedge \alpha_2 \in \Omega^\bullet G$ is homogeneous with weight $w_1+w_2$.
\end{lemma}

We can then generalise Example \ref{ex_thetaiweight}.
\begin{ex}
\label{ex_thetaIweight}
 For any $k=1,\ldots, n,$
   the $k$-form $\theta^I $ given in \eqref{eq_thetaI} is homogeneous with   weight:
 $$
[I]:=\upsilon_{i_1}+\ldots +\upsilon_{i_k}.
$$
In particular, the $n$-form  $\vol $ is homogeneous, and
 its weight is 
 $\upsilon_1+\ldots +\upsilon_n$, also known as the homogeneous dimension of the group $G$. This replace the topological dimension  in the analysis on $G$ by G. Folland and E. Stein  \cite{FollandStein}.
\end{ex}

We denote by $\cW^k$ the set of all the weights occuring in   $\Omega^{k} G$:
$$
 \cW^k 
    :=\{w \ \mbox{is the weight of a form in} \ \Omega^k G \}\ , \quad k\in \bN_0\,.
    $$
    Since $\Omega^0 G = \bR 1$, $\Omega^n G = \bR \vol$ and $\Omega^k G =\{0\}$ for $k>n$, we see:
$$
\cW^n = \{\upsilon_1 +\ldots + \upsilon_n\}\ ,
\ \mbox{ while }\
\cW^k =\{0\} \ \mbox{for}\ k\le0 \ \mbox{or} \ k>n.
$$ 
By Example \ref{ex_thetaIweight},
we know that the forms $\theta^I $'s are homogeneous.
As they form 
the canonical basis of $\Omega^k G$, we obtain the following description in terms of their weights:
$$
    \cW^k 
=\{[I]=\upsilon_{i_1}+\ldots +\upsilon_{i_k}\mid I=(i_1,\ldots,i_k),\ 1\leq i_1<\ldots<i_k \leq n\}\ , 
\quad k=1,\ldots, n\,.
$$
In particular, $\cW^1=\cW(\fg)$ (see \eqref{weights upsilon of g}), and every  $\cW^k$ is a finite (nonempty) subset of $(0,\infty)$ for $k=1,\ldots,n$.   
Furthermore, the space of $k$-forms decomposes through the weights: 
$$
\Omega^k G = \bigoplus_{w\in  \cW^k} \Omega^{k,w} G\ , \qquad k=0,\ldots,n\,,
$$
where we have denoted by
$$
\Omega^{k,w} G := \{ \theta \in \Omega^k G \mid  \delta_r \theta = r^w \theta\ ,\, \forall \,r>0 \}\,,
$$
  the space of  $k$-forms of weight $w\in \cW^k$.

\subsection{On increasing or preserving weights}
\label{subsec_weight}

The formula for $d_0$ in \eqref{eq_d0} together with the dilations $\delta_r$ being automorphisms of the Lie algebra $\fg$ readily imply that  $d_0$ commutes with the dilations 
$$
d_0 (\delta_r \theta) = \delta_r (d_0 \theta) \ , $$
when applied to left-invariant forms
$\theta\in \Hom(\bigwedge^\bullet \fg,\bR)\cong \Omega_L^\bullet G$.
The same formula also holds for any smooth form $\theta\in\Omega^\bullet G$
 since $d_0$ is algebraic
and the $C^\infty(G)$-module $\Omega^\bullet G$ is generated by $\Omega_L^\bullet G $. 
As $d_0$ commutes with the dilations, it preserves the weights of the forms:
$$
d_0(\Omega^{k,w} G)\subset \Omega^{k+1,w} G\ ,\quad k=0,\ldots,n,\, w\in \cW^k.
$$
By construction, this also implies that the maps
$\star$,  $d_0^t$ and $d_0^{-1}$ commute with the action of the dilations.
Hence, they also preserve the weights of the forms.

\begin{lemma}
\label{lem_d-d0}
The operator 
$d-d_0$ strictly increases the weight of a form:
$$
(d^{(k)}-d_0^{(k)})\Omega^{k,w} G\subset \bigoplus_{ w'>w}\Omega^{k+1,w'} G\ ,\quad k=0,\ldots,n\,. 
$$
\end{lemma}

\begin{proof}
By linearity, it is sufficient to consider forms $f\alpha\in\Omega^kG$ with $\alpha\in  \Omega_L^{k,w} G$, $f\in C^\infty(G)$ and $k=0,\ldots,n$. 
By the Leibniz rule, we have
$$
(d-d_0) (f\alpha)= d(f\alpha)-d_0(f\alpha) 
=
(df)\wedge \alpha +f d\alpha - f d_0\alpha
=
(df)\wedge \alpha
,$$
since the  algebraic operator $d_0$ coincides with $d$ on  $\Omega_L^\bullet G$. 
We conclude with Lemma \ref{lem_wedgeweight}, given that $df\in\Omega^1G$ and the set  $\mathcal{W}^1$ of weights occurring in $\Omega^1G$  is a finite subset of $(0,\infty)$.
\end{proof}

Since $\Omega^\bullet G$ contains only a finite number of weights, an operator that strictly increases the weight of a form will necessarily be nilpotent.
Let us formalise this remark in the following statement.
\begin{lemma}
\label{lem_N0}
If $T$ is a linear operator acting on $\Omega^\bullet G$ by strictly increasing the weights, i.e. 	
$$
	T( \Omega^{\bullet,w}G) \subset \bigoplus_{w'>w}  \Omega^{\bullet,w'}G\ ,\ \text{for any weight }w\in \bigcup_{k\in \bN_0}\cW^k\,,
	$$
	then $T$ is nilpotent. 
More precisely, there exists $N_0\in \bN$ independent of $T$ such that $T^{N_0}=0$. 	
\end{lemma}

The integer $N_0$ from Lemma \ref{lem_N0} will of course depend on the weights $\cup_{k\in \bN_0}\cW^k$, and therefore on the Lie structure  of  $\fg$. 
We will keep the same notation for the integer $N_0\in\mathbb N$ from Lemma \ref{lem_N0} throughout the paper. 

It should be noted that the nilpotency of an operator $T\colon\Omega^\bullet G \to\Omega^\bullet G$ that strictly increases the weights will prove useful in showing its invertibility via Neumann series, as pointed out by the following lemma.

\begin{lemma}
\label{lem_DHL4.3}
Let $A_0$ be a linear algebraic operator acting on $\Hom(\bigwedge^\bullet \fg,\bR)$ which naturally extends to $\Omega^\bullet  G$, and let $B$ be a linear differential operator acting on $\Omega^\bullet  G $.
We assume that $A_0$ preserves the weights while $B$ strictly increases the weights.

Then the operator $A:=A_0 -B$
acting on $\Omega^\bullet G$ is 
 invertible, and its inverse $A^{-1}:\Omega^\bullet G \to \Omega^\bullet G$ is a differential operator given by 
\begin{equation}\label{inverse A}
    A^{-1}=  A_0^{-1} \sum_{j=0}^{N_0-1} \big(BA_0^{-1}\big)^j 
	= \bigg(\sum_{j=0}^{N_0-1} \big(A_0^{-1}B\big)^j\bigg)A_0^{-1} ,
\end{equation}
	with $N_0\in\mathbb N$ as in Lemma \ref{lem_N0}.
	\end{lemma}

\begin{proof}
We have
$A= (\id -BA_0^{-1}) A_0 $.
The differential operator
$BA_0^{-1}$  acts on $\Omega^\bullet G$ and strictly increases the weights.
By Lemma \ref{lem_N0},  $(BA_0^{-1})^{N_0} =0$.
Moreover $\id -BA_0^{-1}$ is invertible, 
and its inverse is given by $\sum_{j=0}^{N_0-1} (BA_0^{-1})^j$. Since $A^{-1}=A_0^{-1}(\id-BA_0^{-1})^{-1}$, we obtain the first equality of \eqref{inverse A}. To obtain the second one, it is sufficient to apply the same reasoning to the expression $A= A_0(\id -A_0^{-1}B)  $.
\end{proof}

\subsection{Rumin's construction of $d_c$}

Using
the previous considerations on the invertibility of differential operators that strictly increase the weights of forms, 
Rumin constructs the projection $\Pi$ in the following way.

	\begin{lemma}
	\label{lem_Pi} Let us consider the differential operator $b:=d_0^{-1}d_0 - d_0^{-1} d=- d_0^{-1}(d-d_0)$ acting on $\Omega^\bullet G$.
\begin{enumerate}
\item The map $b:\Omega^\bullet G \to \Omega^\bullet G $ is  nilpotent, and therefore $\id-b$ is invertible and $(\id-b)^{-1}$ is a well-defined differential operator acting on $\Omega^\bullet G$.
\item 
 The differential operator 
		 \begin{align}\label{def proj Pi}
	\Pi:=(\id-b)^{-1} d_0^{-1}d+d(\id-b)^{-1} d_0^{-1} 
	\end{align}
	 is the projection of $\Omega^\bullet G$ onto
	$$
	 F:=\IM d_0^{-1} +  \IM dd_0^{-1}=\IM d_0^{t} +  \IM dd_0^{t}\,,
	 $$
	 along 
$$
 E:=\ker  d_0^{-1} \cap  \ker (d_0^{-1}d)
= \ker  d_0^{t} \cap  \ker (d_0^{t}d)\,.
$$
\end{enumerate}
\end{lemma}

\begin{proof}
By Lemma \ref{lem_d-d0}, $b$ strictly increases weights,
so it is nilpotent, and in particular, by Lemma \ref{lem_N0}, there exists $N_0\in\mathbb N$ such that $b^{N_0}=0$.
We can now apply Lemma \ref{lem_DHL4.3} with $B=b$ and $
A_0=\id\,$: the operator $\id-b$ is invertible, and its inverse is given by the differential operator
\begin{align}\label{formula (1-b) inverse}
  (\id-b)^{-1} = \sum_{j=0}^{N_0-1} b^j = \id +b + \ldots + b^{N_0-1}=\sum_{j=0}^{N_0-1}(-1)^j\big[d_0^{-1}(d-d_0)\big]^j\,.  
\end{align}
This implies that $\Pi$ is a well-defined differential operator on $\Omega^\bullet G$, and concludes the proof of part (1).

By \eqref{formula (1-b) inverse} and \eqref{def proj Pi}, we have that the image of $\Pi$ is included in the subspace $F\subset\Omega^\bullet G$ defined in the statement.
Moreover, the properties of $d,d_0$ and $d_0^{-1}$  imply
$$
\Pi d_0^{-1} = d_0^{-1}
\quad\mbox{and}\quad
\Pi dd_0^{-1}=  dd_0^{-1}\ ,
$$
and so $\Pi=\id$ on $F$ and $\Pi^2=\Pi$, i.e. $\Pi$ is a projection onto $\IM \Pi = F$ along $\ker \Pi$.

Clearly, 
$\ker \Pi$ contains the subspace $E$ defined in the statement. 
The properties of $d,d_0$ and $d_0^{-1}$ also imply
$$
d_0^{-1} \Pi 
= d_0^{-1} \ \text{ and }
\ 
d_0^{-1}d \Pi 
=d_0^{-1}d\,,
$$
so $\ker\Pi =E$. 
This shows part (2) and concludes the proof.
\end{proof}

Following Rumin's notation, we  denote the two projections onto $F$ and $E$ respectively as
$$
\Pi_F:=\Pi 
\qquad\mbox{and}\qquad 
\Pi_E:=\id -\Pi\, ,
$$
where $\Pi$ is defined in \eqref{def proj Pi}.
It is easy to check that they enjoy the following properties:
\begin{lemma}
	\label{lem_Pis}
	We have
\begin{itemize}
\item[1.]	
$d_0^{-1}\Pi_E = \Pi_E d_0^{-1}=0$;
\item[2.]
$d\Pi_F = \Pi_F d$, and $d\Pi_E = \Pi_E d$;
\item[3.]
$\Pi_E\pr_{E_0^\perp}\Pi_E =0$;
\item[4.]
 on $\Omega^k G$, we have
$\Pi_F^t=(-1)^{k(n-k)}\star \Pi_F \star $ and $\Pi_E^t=(-1)^{k(n-k)}\star \Pi_E \star $.
\end{itemize}
\end{lemma}
Note that 4. implies $F^\perp = \star F \star$ and $E^\perp = \star E \star$. 

\begin{proof}
	The equality in 1. follows from the fact that $d_0^{-1}\Pi=\Pi d_0^{-1}=d_0^{-1}$.
	The equalities in 2. can be easily checked by carrying out the explicit computations that show $d\Pi=d(\id-b)^{-1}d_0^{-1}d=\Pi d$.
	For the equality in 3. it is sufficient to apply property 1., since $\pr_{E_0^\perp}=\id-\Pi_0 = d_0^{-1}d_0 +d_0d_0^{-1}$ by Proposition \ref{prop_prE0}, and so
$$
\Pi_E\pr_{E_0^\perp}\Pi_E = \Pi_E(d_0^{-1}d_0 +d_0d_0^{-1})\Pi_E =0.
$$

Finally, let us observe that the formal transpose of $d-d_0$ and $d_0$ coincide with $\star (d-d_0)\star$ and 
$\star d_0\star$ up to a sign on $\Omega^\bullet G$, and a similar property is true for $b^j$, $j=0,1,\ldots,N_0-1$, and hence for $(\id-b)^{-1}d_0^{-1}$. Therefore, $\Pi^t$ coincides with $\star \Pi \star $ up to a sign on $\Omega^k G$. 
We compute 
$$
(\star \Pi \star)^2
=\star \Pi \star\star \Pi \star
=(-1)^{k(n-k)}\star \Pi \star
=(-1)^{k(n-k)}\star \Pi \star
$$
As $\Pi^t$ is the  projection $\Pi_F^t$, 
we obtain the formulae in 4.   
\end{proof}

Lemma \ref{lem_Pis} implies readily the following result.

\begin{theorem}[M. Rumin]
\label{thm_Rumin}
   The de Rham complex $(\Omega^\bullet, d)$ splits into two subcomplexes $(E^\bullet,d)$ and $(F^\bullet,d)$. 
Moreover,   the differential operator acting on $\Omega^\bullet G$ and defined by 
$$
d_c :=\Pi_0\, d\, \Pi_E\, \Pi_0 
$$
satisfies 
$$
d_c^2=0\quad ,\quad d_c (\Omega^k G)\subset \Omega^{k+1} G\quad 
\text{ and }  \quad 
d_c^{(k,t)} = (-1)^{kn+1} \star\, d_c \, \star \
\mbox{on}\ 
\Omega^{k+1}\, .
$$ 

The  complex  $(E_0^\bullet,d_c)$ computes the de Rham cohomology of the underlying homogeneous group $G$ and it is known as the Rumin complex.
\end{theorem}

\section{An alternative construction of the Rumin complex}
\label{sec_altconst}

In this section, we present an alternative construction of the Rumin complex $(E_0^\bullet,d_c)$ with ideas from parabolic geometry. 
Importantly,  this new construction will not rely on the partial inverse $d_0^{-1}$ of the algebraic map $d_0$.

\subsection{The operators $\Box$, $P$ and $L$}
\label{subsec_BoxPL}

Our box operator $\Box$ is not the usual de Rham Laplace operator, but instead the differential operator of order 1 (in the Euclidean sense) given by
$$
\Box := d d_0^t  +d_0^t d\,.
$$
In general, it may not have any particular properties such as symmetry or even commutation with its formal transpose.
However, just like $\Box_0$, $\Box$ maps $\Omega^k G$ to itself for any degree $k$ with 
$$
\Box=\sum_{k=0}^n \Box^{(k)}\ , \quad 
\Box^{(k)} :=  d^{(k-1)} d_0^{(k-1,t)} +d_0^{(k,t)} d^{(k)}:
\Omega^k G \longrightarrow \Omega^k G\,.
$$
Furthermore,  we have
$$
\Box = \Box_0 \ + \   (d-d_0)d_0^t +d_0^t (d-d_0)\ ,
\quad\mbox{since}\quad
\Box_0 = d_0 d_0^t +d_0^t d_0\,.
$$
The operators $\Box_0$ and $d_0$ are algebraic and preserve the weights of the forms, while the differential operator $(d-d_0)$ strictly increases the weights, see Section \ref{subsec_weight}.
Therefore, $\Box-\Box_0$ also strictly increases the weights. 

The operator $\Pi_0$ is the spectral projection of the symmetric linear map $\Box_0$ onto its kernel, with both  operators viewed as acting on $\Hom(\bigwedge^\bullet \fg,\mathbb R)$, and then extended as algebraic operators acting on $\Omega^\bullet G$. 
As we do not have a spectral decomposition for $\Box$, we cannot define analogously an operator $P$ as a spectral projector in the same way.
However, we notice that 
\begin{equation}
	\label{eq_Pi0_int}
	\Pi_0 = \frac1{2\pi i}\oint_{|z|=\varepsilon} (z-\Box_0)^{-1} dz ,
\end{equation}
where $|z|=\varepsilon$ denotes the circle about 0 with radius  $\eps>0$ small enough, see point 2. of Lemma \ref{lem_matrix}.
This expression is independent of $\eps\in (0,\lambda_1)$ where $\lambda_1>0$ is the first non-zero eigenvalue of $\Box_0$.
We then define the operator $P$  in a similar manner but for $\Box$:
we set, at least formally, 
\begin{equation}
	\label{eq_defP}
P: = \frac1{2\pi i}\oint_{|z|=\eps} (z-\Box)^{-1} dz \ ,
\end{equation}
with $\varepsilon>0$  as in the formula \eqref	{eq_Pi0_int} for $\Pi_0$. 

Let us explain why the definition of $P$ via \eqref{eq_defP} makes sense. 
We write  for any $z\in \bC$, 
$$
z-\Box = z-\Box_0 - B, 
\quad\mbox{where}\quad 
B:=\Box - \Box_0=(d-d_0)d_0^t +d_0^t (d-d_0)\,.
$$
We observe that, for $z\in \bC \setminus \spec (\Box_0)$, 
$A_0=z-\Box_0$ is an invertible algebraic operator that preserves the weights while, as already mentioned, the differential operator $B$ strictly increases the weights.
By Lemma \ref{lem_DHL4.3}, $z-\Box$ is invertible, and its inverse is given by the finite sum 
 	\begin{equation}
 		\label{eq_resolventBox}
 		(z-\Box)^{-1}= (z-\Box_0)^{-1}  \sum_{j=0}^{N_0-1} \big(B(z-\Box_0)^{-1}\big)^j\,.
 	\end{equation}
In particular, $(z-\Box)^{-1}$ is a differential operator on $\Omega^\bullet G$
whose coefficients are meromorphic functions of $z$ with poles lying in  $\spec(\Box_0)$. 
Consequently, \eqref{eq_defP} makes sense and defines a differential operator $P$ on $\Omega^\bullet G$
independently of $\eps\in (0,\lambda_1)$. 

The operator $P$ enjoys  the following important properties.
\begin{proposition}
\label{prop_P}
The formula in \eqref{eq_defP} defines a differential operator $P$ acting on $\Omega^\bullet G$ that satisfies the following properties:
\begin{itemize}
	\item[i.] $P$ is a projection, i.e. $P^2=P$. It commutes with $\Box$, $d$ and $d_0^t$, and preserves the degrees of the forms: $P(\Omega^k G) \subset \Omega^k G$;
	\item[ii.] $P-\Pi_0$ strictly increases weights and $\Pi_0 P \Pi_0 =\Pi_0$;
\item[iii.] 
the operator $\Box$ acts on  $\IM P$ and on $\ker P$.
Moreover, it is nilpotent on $\IM P$ with $\Box^{N_0} P=0$, 
and it is  invertible on $\ker P$.
Consequently, $P$ is the projection onto $\ker \Box^{N_0}$ along $\IM \Box^{N_0}$. 
\end{itemize}	
\end{proposition}

Property iii. above should be compared with the analogous property of the algebraic operator $\Box_0$, i.e. $\Box_0$ is zero on $\IM \Pi_0= \ker \Box_0$, it is invertible on $\ker \Pi_0 = \IM \Box_0$, and hence $\Pi_0$ is the projection onto $\ker \Box_0$ along $\IM \Box_0$.

\begin{proof}[Proof of i. and ii.]
By the Cauchy integral formula, when $\eps_2>\eps_1$, we have
\begin{align*}
	P^2 &=  \frac1{(2\pi i)^2}\oint_{|z_1|=\eps_1}\oint_{|z_2|=\eps_2} (z_1-\Box)^{-1}(z_2-\Box)^{-1} dz_1 dz_2\\
	&=  \frac1{(2\pi i)^2}\oint_{|z_1|=\eps_1}\oint_{|z_2|=\eps_2} \frac 1{z_2-z_1} \big( (z_1-\Box)^{-1}-(z_2-\Box)^{-1} \big)dz_1 dz_2\\
	&=  \frac1{(2\pi i)^2}\oint_{|z_1|=\eps_1}\oint_{|z_2|=\eps_2} \frac 1{z_2-z_1}  dz_2 \ (z_1-\Box)^{-1}dz_1
	= \frac1{2\pi i}\oint_{|z_1|=\eps_1}(z_1-\Box)^{-1}dz_1 \,.
\end{align*}
Since $\Box$ commutes with its resolvant $(z-\Box)^{-1}$, 
it also commutes with $P$.
As $\Box$ commutes with $d$ and $d_0^t$ and  preserves the degrees of the forms, so does its resolvant, and then $P$. 
This shows i. 

 The formulae \eqref{eq_defP}, \eqref{eq_Pi0_int}, \eqref{eq_resolventBox} for $P$, $\Pi_0$ and $(z-\Box)^{-1}$ yield:
\begin{align*}
	P&=\frac1{2\pi i}\oint_{|z|=\epsilon}(z-\Box_0)^{-1}  \sum_{j=0}^{N_0-1} \big(B(z-\Box_0)^{-1}\big)^j
 dz \\
 &=\Pi_0
 +\sum_{j=1}^{N_0-1}\frac1{2\pi i}\oint_{|z|=\epsilon}(z-\Box_0)^{-1}
   \big(B(z-\Box_0)^{-1}\big)^j
 dz\,,
\end{align*}
which shows that $P-\Pi_0$ strictly increases the weight since $B $ does.

By the spectral decomposition of $\Box_0$, we have
$$
\Pi_0 (z-\Box_0)^{-1}
   B(z-\Box_0)^{-1} \Pi_0
   = \frac 1{z^2} \Pi_0 B\Pi_0\,,
$$
and more generally for $j\geq 1$
$$
\Pi_0 (z-\Box_0)^{-1}
   \big(B(z-\Box_0)^{-1}\big)^j
   \Pi_0 = \frac 1{z^2}
   \Pi_0 
   \big(B(z-\Box_0)^{-1}\big)^{j-1}
   \Pi_0\,.
 $$  
Hence, the residue of  
$\Pi_0 (z-\Box_0)^{-1}
   \big(B(z-\Box_0)^{-1}\big)^j \Pi_0 $ at $z=0$ is zero when $j>0$.
   Consequently, $\Pi_0 P \Pi_0 = \Pi_0^3 = \Pi_0$. This concludes the proof of Part  ii.
\end{proof}

The proof of property iii. will be given at the end of this section. It will use the operator $L$  defined and studied in the next statement:

\begin{proposition}
\label{prop_L}
The differential operator $L$ acting on $\Omega^\bullet G$ and defined as
\begin{align}\label{def operator L}
L:=P\Pi_0 +(\id -P)(\id-\Pi_0)    
\end{align}
preserves the degrees of the forms.
It is invertible and its inverse is a differential operator acting on $\Omega^\bullet G$.
We have 
$$
PL=P\Pi_0 = L\Pi_0\quad,\quad\text{ so }\quad
\Pi_0 = L^{-1} P L\,.
$$
Moreover, $L-\id$, $L^{-1}-\id$ and $L^{-1} \Box L - \Box_0$ strictly increase weights, and $L^{-1} \Box L$ commutes with $\Pi_0$, i.e.
$\Pi_0 L^{-1} \Box L = L^{-1} \Box L\Pi_0$. 
\end{proposition}

\begin{proof}[Proof of Proposition \ref{prop_L}]
The differential operator $P$ and the algebraic operator $\Pi_0$ preserve the degrees of the forms, thus so does $L$.
As $P$ and $\Pi_0$ are projections,  we easily check that 
\begin{equation}
	\label{eq_PlPi0} 
	PL = P\Pi_0\qquad \mbox{and}\qquad
L\Pi_0= P\Pi_0\,,
\end{equation}
and that
$$
L - I = (P-\Pi_0)\Pi_0 + (\Pi_0-P)(I-\Pi_0)
 \,,
$$
which shows that
$L-\id$ strictly increases weights 
by property ii. of Proposition \ref{prop_P}.  Hence  $L$ is invertible by Lemma \ref{lem_DHL4.3}, 
with 
$$
L^{-1} = \sum_{j=0}^{N_0-1} (\id -L)^j.
$$
Therefore, $L-\id$, $L^{-1}-\id$ and $\Box-\Box_0$ strictly increase weights, while $\id$ and $\Box_0$ keep the weights constant. 
This readily implies that 
$L^{-1} \Box L - \Box_0$ strictly increases the weights. 
Finally, from \eqref{eq_PlPi0} and $P\Box=\Box P$, we have that the two operators $L^{-1} \Box L$ and $\Pi_0$ commute. 
\end{proof}

We can now finish the proof of Proposition \ref{prop_P}.
\begin{proof}[Proof of iii. of Proposition \ref{prop_P}]
As $P$ and $\Box$ commute, 
we already know that $\Box$ acts on $\IM P$ and $\ker P$, 
but we can be more precise. 
As $L^{-1} \Box L$ and $\Pi_0$ commute, 
$L^{-1} \Box L$ acts on $\IM \Pi_0 = \ker \Box_0$ and on $\ker \Pi_0 = \IM \Box_0$. 
As $L^{-1} \Box L - \Box_0$ strictly increases the weights,
we have 
by Lemma \ref{lem_N0}
\begin{equation}
	\label{eq_LBoxN0}
	(L^{-1} \Box L - \Box_0)^{N_0}=0\,.
\end{equation}
Considering \eqref{eq_LBoxN0} on $\IM \Pi_0 = \ker \Box_0$, we obtain 
$0= (L^{-1} \Box L)^{N_0}= L^{-1} \Box^{N_0} L$, 
hence $\Box^{N_0}=0$ on $L \IM \Pi_0 L^{-1} = \IM L \Pi_0 L^{-1} = \IM P$. 

We observed previously that 
the operators  $L^{-1} \Box L$, $\Box_0$ and $L^{-1} \Box L-\Box_0$ act on $\ker\Pi_0 = \IM \Box_0$.
We  denote by $\tilde A$, $\tilde A_0$ and $\tilde B$ respectively the operators  $L^{-1} \Box L$, $\Box_0$ and $L^{-1} \Box L-\Box_0$ restricted to and acting on $\ker\Pi_0$.
In particular, $\tilde A_0$ is invertible, and by \eqref{eq_LBoxN0} we have $\tilde B^{N_0} = 0$.
Proceeding as in the proof of Lemma \ref{lem_DHL4.3}, we see that $\tilde A$ is also invertible.
This implies that  $L$ is invertible on $L \ker \Pi_0 L^{-1} = \ker L \Pi_0 L^{-1} = \ker P$.
\end{proof}

Let us examine the effect of $\star$-conjugation  on the objects constructed in this section.
\begin{lemma}
\label{lem_star_BoxPL} 
We have
$$
\star \Box^{(k)} \star = (-1)^{k(n-k)}  \Box^{(n-k,t)}\quad ,
\quad
\star P^{(k)} \star = (-1)^{k(n-k)} P^{(n-k,t)}\ ,
$$
and
$$
\star L^{(k)} \star = (-1)^{k(n-k)} L_1^{(n-k)}\ ,
\quad\mbox{where}\quad 
L_1:=P^t \Pi_0 +(\id -P^t)(\id-\Pi_0)\,.    
$$
The linear map $L_1:\Omega^\bullet G\to \Omega^\bullet G$ is a differential operator that preserves the degree of forms. It is invertible and its inverse is a differential operator.
\end{lemma}
\begin{proof}
From the properties of $\star$, $d$ and $d_0$
(see Sections \ref{subsec_extOmega} and \ref{subsubsec_star0}), 
we compute
\begin{align*}
	\star \Box^{(k)} \star &= (-1)^{(k-1)(n-k+1)}
	\star d^{(k-1)} \star \star d_0^{(k-1,t)} \star 
	+(-1)^{(k+1)(n-k-1)} \star d_0^{(k,t)} \star \star d^{(k)}\star 
	\\
	&=  (-1)^{k(n-k)} d^{(k-1,t)} d_0^{(k-1)} +(-1)^{k(n-k)}  d_0^{(k)} d^{(k,t)}=  (-1)^{k(n-k)}  \Box^{(n-k,t)}\,.
\end{align*}
The $\star$-conjugation for $P$ follows from this.
As $\Pi_0$ commutes with $\star$, we obtain the formula that relates $L$ to $L_1$. The properties of the operator $L_1\colon\Omega^\bullet G\to\Omega^\bullet G$ can be shown by following the same steps as in the proof of Proposition \ref {prop_L} for the same properties for $L$ .
\end{proof}

\subsection{Decoupling of the complex chain $L^{-1}dL$}
As  the de Rham differential	 $d$ commutes with  $P$ (see point i. of Proposition \ref{prop_P}), $L^{-1} dL $ commutes with $L^{-1}PL$, the latter being equal to $\Pi_0$ by Proposition \ref{prop_L}. 
The commutation of  $\Pi_0$  with $L^{-1} dL$ allows us to decompose the differential chain $L^{-1} d L$ as
$$
L^{-1} d L \ = \ L^{-1} d L \Pi_0  + L^{-1} d L (\id -\Pi_0)
\ = \ D+C\,,
$$
where the differential operator 
\begin{align}\label{def operator D}
    D:= L^{-1} d L \Pi_0 \ =\ \Pi_0 L^{-1} d L \Pi_0
\end{align}
acts on $\IM \Pi_0$, 
while the differential operator 
\begin{align}\label{def operator C}
    C:= L^{-1} d L (\id -\Pi_0)\ =\ (\id-\Pi_0) L^{-1} d L (\id -\Pi_0)
\end{align}
acts on $\ker \Pi_0$. 

In other words, 
the chain complex $L^{-1}dL$ decomposes into the direct sum of two chain complexes: one where the operator $D$ acts on $E_0=\IM \Pi_0$, and the other where $C$ acts on $F_0=\ker\Pi_0$. In fact,
$$
(L^{-1}dL)^2 =0
\quad\text{ and }\quad
L^{-1}dL = D+C
\quad\mbox{with}\quad
D^2=C^2=0\,. 
$$ 

\subsubsection{The chain complex $(F_0^\bullet,C)$} As already pointed out, the operator $C$ defined in \eqref{def operator C} acts on 
$\ker \Pi_0 = F_0$.
It is not difficult to construct an explicit invertible operator conjugating the maps $C$ and $d_0$ on $F_0$, and this implies that the cohomologies of the complexes $(F_0^\bullet,C)$ and  $(F_0^\bullet,d_0)$ are isomorphic and therefore both trivial.
\begin{proposition}
\label{prop_GconjBd0}
\begin{enumerate}
\item The maps $C$ and $d_0$ are conjugated on $F_0$:
$$
Cg=gd_0 \quad\mbox{on}\  F_0\,, 
$$
 where $g:F_0^\bullet \to F_0^\bullet$ is the invertible operator defined as
$$
g:= C(\id-\Pi_0)d_0^t \Box_0^{-1} +(\id-\Pi_0) d_0^t \Box_0^{-1} d_0 \quad\mbox{on}\quad F_0=\ker \Pi_0\,.
$$
    \item The complex $(F_0^\bullet,C)$ has trivial cohomology, i.e. for any $k=0,\ldots,n$
$$
H^k (F_0,C) = \frac{\ker C\colon F_0^k\to F_0^{k+1}}{ \IM C\colon F_0^{k-1}\to F_0^k} =\{0\}\,.
$$
\end{enumerate}
\end{proposition}

\begin{proof}[Proof of Proposition \ref{prop_GconjBd0}]
The operator $g$ given in the statement is well defined as an operator acting on $F_0$. 
Since $L-\id$, $L^{-1}-\id$ and $d-d_0$ strictly increase the weights, while $\Pi_0$ and $\Box_0$ keep the weights constant,
 the operator 
 $$
 C-(\id -\Pi_0)d_0 (\id-\Pi_0) 
 = (\id -\Pi_0) L^{-1} d L(\id-\Pi_0)
 -(\id -\Pi_0)d_0 (\id-\Pi_0) 
 $$
 strictly increases the weights, and so does the  operator
$$
g - (\id -\Pi_0)d_0 (\id-\Pi_0) (\id-\Pi_0)d_0^t \Box_0^{-1} +(\id-\Pi_0) d_0^t \Box_0^{-1} d_0
\quad \mbox{on} \ \ker \Pi_0\,.
$$
By Proposition \ref{prop_prE0}
and the fact that $d_0$ commutes with $\Box_0$ and therefore also with its inverse,
the latter operator simplifies as 
$$ 
g -(\id -\Pi_0) (d_0 d_0^t +  d_0^td_0) \Box_0^{-1}
=g - (\id -\Pi_0).
$$
Hence, when restricted to $\ker\Pi_0$, $g - \id$ strictly increases weights,
and so if we proceed as in the proof of Lemma \ref{lem_DHL4.3}, we have that $g$ is invertible on $F_0=\ker \Pi_0$. 

To check that $Cg =gd_0$ on $F_0 = \IM d_0 +\IM d_0^t$, 
 it suffices to check that the equality holds both on $\IM d_0$ and on $\IM d_0^t$. 
As $d_0^2=0$ and $C^2=0$, we see that both $Cg$ and $gd_0$ vanish on $\IM d_0$, while both 
$gd_0 d_0^t $ and $ Cg d_0^t $ are equal to $C (\id-\Pi_0) d_0^t \Box_0^{-1}  d_0d_0^t$.
This shows Part (1).

Finally, by Part (1) we have that $g$ induces a linear isomorphism
between the cohomology groups of the two complexes $(F_0^\bullet,C)$ and  $(F_0^\bullet,d_0)$ given by
$$
g_H \colon\left\{\begin{array}{rcl}
     H^\bullet(F_0,d_0)&\longrightarrow& H^\bullet(F_0,C)\\
\ [\alpha]  & \longmapsto &  [g \alpha] 
\end{array}\right. ,
$$
with 
inverse given by $[\beta] \mapsto [g^{-1}\beta]$. Part (2) then follows directly from Lemma \ref{lem_Hd0F0}, where it was shown that the cohomology of the complex $(F_0^\bullet,d_0)$ is trivial.
\end{proof}

\subsubsection{The chain complex $(E_0^\bullet,D)$}

The differential operator $D$ defined in \eqref{def operator D} acts on 
$E_0=\IM \Pi_0$. In this subsection, we show that the complex $(E_0^\bullet,D)$ is conjugated to the de Rham complex $(\Omega^\bullet G,d)$ via the operator $L$ introduced in \eqref{def operator L}, i.e.
$$
L D = d  L
\quad\mbox{on} \ E_0\,.
$$	
Indeed, this follows easily from 
the equalities $L\Pi_0L^{-1} = P$ and $PL \Pi_0=L\Pi_0$ proved in Proposition \ref{prop_L}, and the fact that $d$ commutes with $P$ by Proposition \ref{prop_P}.
Hence, 
the differential operator $L\Pi_0$ is a chain map between $(\Omega^\bullet G, d)$ and $(E_0^\bullet,D)$, i.e.
$(L\Pi_0) \, D  =d\, (L\Pi_0)$, and it induces the linear map 
$$
(L\Pi_0)_H\colon \left\{\begin{array}{rcl}     H^\bullet(E_0,D)&\longrightarrow & H^\bullet(\Omega^\bullet G,d)\\
     \ [\alpha]   &\longmapsto &  [L\Pi_0\alpha]
     \end{array}\right. ,
     $$
between the cohomology groups of $(E_0^\bullet,D)$ and $(\Omega^\bullet G,d)$.

This map $(L\Pi_0)_H$ is unlikely to be  a linear isomorphism in general. 
However, the next statement implies that the map induced by $\Pi_0 L^{-1}$ and given by
$$
(\Pi_0L^{-1})_H\colon \left\{
\begin{array}{rcl}
     H^\bullet(\Omega^\bullet G,d)& \longrightarrow& H^\bullet(E_0,D)\\
 \ [\alpha]&\longmapsto & [\Pi_0L^{-1}\alpha]
\end{array}\right. ,
$$
provides such an isomorphism, with inverse given by $[\beta]\mapsto [L\beta]$.

\begin{proposition}
The differential operator $\Pi_0 L^{-1}$ 
is a chain map between $(\Omega^\bullet G,d)$ and $(E_0^\bullet,D)$, that is $(\Pi_0 L^{-1})d = D(\Pi_0 L^{-1})$.
This chain map is homotopically invertible, with homotopic inverse given by $L$
since we have
$$
(\Pi_0 L^{-1}) L =\Pi_0,
\quad\mbox{and}\quad
\id - L(\Pi_0 L^{-1}) = d h +hd\,,
$$
where $h$ is the differential operator acting on $\Omega^\bullet G$ and defined as
$$
h:=L\, g (\id -\Pi_0)\, d_0^t\Box_0^{-1} \, g^{-1}(\id -\Pi_0)\, L^{-1}\,.
$$
\end{proposition}

\begin{proof}
	By construction, we have
	$d = L(D+C)L^{-1}$, with $D=\Pi_0 D\Pi_0$ and $C=(\id-\Pi_0) C (\id-\Pi_0)$,
	so 
	$\Pi_0 L^{-1} d = \Pi_0 D \Pi_0 L^{-1}= D \Pi_0 L^{-1}$.
	
	By applying the definitions of $h$ and $C$ together with Proposition \ref{prop_GconjBd0}, we have
	\begin{align*}
	L^{-1} (dh+hd) L&=
	C gd_0^t \Box_0^{-1} \, g^{-1}(\id -\Pi_0) 
	+ g (\id -\Pi_0)\, d_0^t\Box_0^{-1} \, g^{-1} C\\
	&= g\left( d_0 d_0^t \Box_0^{-1} + d_0^t \Box_0^{-1}d_0\right) g^{-1} (\id-\Pi_0)\,.
	\end{align*}
	As $\id_{\ker \Pi_0}= d_0 d_0^t \Box_0^{-1} + d_0^t \Box_0^{-1}d_0 $, we have obtained
	$L^{-1} (dh+hd) L=\id-\Pi_0$. The statement follows. 
\end{proof}

\begin{corollary}
The cohomology of $(E_0^\bullet,D)$ is linearly isomorphic to the de Rham cohomology of the manifold $G\cong\bR^n$:
$$
H^0 (E_0,D) \cong \bR\quad ,
\quad \mbox{and for}\ k=1,\ldots,n, \quad 
H^k(E_0,D)=\{0\}.
$$
\end{corollary}

We conclude this section by examining the $\star$-conjugation of $D$.
\begin{lemma}
\label{lem_starD}
We have
$$
D^{(k,t)} = (-1)^{kn+1} \star\, D^{(n-k-1)} \, \star \ .
$$
\end{lemma}
\begin{proof}
By Lemma \ref{lem_star_BoxPL} and the commutation of $\star$ with $\Pi_0$, we have on $\Omega^{k+1} G$
$$
\star D^{(n-k-1)} \star  
= (-1)^{kn+1}\Pi_0^{(k)} L_1^{(k,-1)} d^{(k,t)} L_1^{(k+1)} \Pi_0^{(k+1)}\,.
$$
We would like to compare this with 
$D^t = \Pi_0 L^t d^t L^{-t} \Pi_0$ with $L^{-t} = (L^{-1})^t = (L^t)^{-1}$.
We observe that since $P$ and $\Pi_0=\Pi_0^t$ are projections,  we have
\begin{align*}
L^t L_1 \Pi_0 &= L^t P^t \Pi_0  = \Pi_0 P^t P^t \Pi_0 
= (\Pi_0 P \Pi_0)^t=(\Pi_0)^t =\Pi_0\,,
\end{align*}
by point ii. of Proposition \ref{prop_P}, and similarly
$$
\Pi_0 L^t L_1 = \Pi_0  P^t L_1 =  \Pi_0 P^t P^t \Pi_0  = \Pi_0\,.
$$
Therefore, $L_1 \Pi_0 = L^{-t}\Pi_0$ and $\Pi_0 L^t = L_1^{-1}\Pi_0$, and so the formula follows.
\end{proof}

\subsection{The two constructions coincide}

In Section \ref{sec_licohomG}, 
we have established 
that the  orthogonal  projection $\Pi_0$ onto $E_0= \ker d_0^t \cap \ker d_0$ is the same in Rumin's construction and as  the spectral projection of $\Box_0$ for the 0-eigenvalue.
Here, we are going to show that Rumin's complex $(E_0^\bullet,d_c)$ constructed by using $d_0^{-1}, \Pi_E$ etc. and the complex $(E_0^\bullet,D)$ constructed via $\Box, P,L$ etc. coincide.

\begin{lemma}
\label{lem_PiE=P}
The projections $\Pi_E=\id-\Pi$ and $P$ coincide.
\end{lemma}

\begin{proof}
As $\Pi_0 d_0^t=0$ and $P$ commutes with $d_0^t$ 
(see Proposition \ref{prop_P}), we have
$$
(P-\Pi_0)d_0^t =Pd_0^t = d_0^t P\,,
$$
and so recursively for any $k\ge 1$
$$
(P-\Pi_0)^{(k)}d_0^{(k-1,t)} =P^{(k)}d_0^{(k-1,t)}= d_0^{(k-1,t)} P^{(k)}\,.
$$
Let us stress that for $k\ge N_0$, the left-hand side vanishes by Proposition \ref{prop_P}.
As $P$ is a projection, we have obtained 
$0 =Pd_0^t =d_0^t P$. Moreover, since $P$ commutes with $d$ (see Proposition \ref{prop_P}), we also have $Pdd_0^t = dP d_0^t =0$, which implies the inclusion $\ker P \supset \IM d_0^t + \IM dd_0^t =F$.

By point iii. of Proposition \ref{prop_P}, $P$ is the projection onto $\ker \Box^{N_0}$ along $\IM \Box^{N_0}$.
Since $d^2=0$ and $(d_0^t)^2=0$, we compute easily 
$$
\Box^{N_0} = (d d_0^t)^{N_0} +(d_0^t d)^{N_0},
$$
so we have 
$$
\IM P =\ker \Box^{N_0}  \supset \ker d_0^t \cap \ker (d_0^t d) = E\quad 
 \mbox{and}\quad
\ker P =\IM \Box^{N_0}  \subset \IM d_0^t + \IM (d d_0^t)= F\,.
$$
These last inclusions then imply that $\ker P= F$, but also $\IM P =E$, and so $P=\Pi_E$. 
\end{proof}

Clearly,  both $d_c$ and $D$ increase the degrees of the forms by 1, and we have already obtained that they behave in the same way under $\star$-conjugation (see Lemma \ref{lem_starD} and Theorem \ref{thm_Rumin}).
We will now prove that these  maps coincide:
\begin{theorem}
	The  maps $d_c$ and $D$ coincide: 
$d_c = D$.
\end{theorem}

\begin{proof}
By Proposition \ref{prop_L}, the construction of $L$ implies $\Pi_0 L^{-1}P\Pi_0 = \Pi_0$, and since $P (1-\Pi_0)P=0$ by point 3. of Lemma \ref{lem_Pis}, we have
$$
\Pi_0 L^{-1} = \Pi_0 L^{-1} P^2
= \Pi_0 L^{-1} P (\Pi_0 + (\id-\Pi_0)) P
= \Pi_0 L^{-1} P \Pi_0  P = \Pi_0 P\,.
$$
This equality, together with the fact that $d$ commutes with $P$
by point 2. of Lemma \ref{lem_Pis}, and that $PL\Pi_0 = P\Pi_0$ by Proposition \ref{prop_L}, implies
$$
D= \Pi_0 L^{-1} d L \Pi_0 
= \Pi_0 P d L \Pi_0
=\Pi_0 d PL\Pi_0 
=\Pi_0 d P\Pi_0\,,
$$
which coincides with $d_c$
since $P=\Pi_E$ by Lemma \ref{lem_PiE=P}. 	
\end{proof}

\section{The subcomplex $(E_0^\bullet,D)$ on an isomorphic family of Engel groups}

The aim of this section is to present the explicit computations behind the construction of the Rumin complex that was described in Section \ref{sec_altconst} for the specific case of the Engel group. In order to do this, we will give concrete matrix  representations of the different maps that are involved in the construction of the differential operator $D$.

This work will be done for a 1-parameter family $G_t$ of stratifiable groups which are all isomorphic to the Engel group for each $t>0$. Surprisingly, we will see how some parts of the construction, for instance the spectral decomposition of $\Box_0$,
 heavily rely on the parameter $t>0$.
This highlights the fact that, even though the complexes $(E_0^\bullet,D)$ and $(F_0^\bullet,C)$ are essentially unique up to Lie algebra isomorphisms, the intermediate steps that lead to the construction of these subcomplexes depend on the parameter $t$, and hence on the particular realisation of the group.

The case of $G_t$ for $t=1$ was given in  \cite[Appendix A]{DaveHaller}.

\subsection{The Engel groups $G_t$ and the canonical basis on the exterior algebra}
\label{subsec_GtthetaI}

Here $t>0$ is a fixed positive number.
We denote by $\fg_t$ the four dimensional 
Lie algebra $\fg_t$ equipped with the stratification
$$
\fg_t = \fg_1 \oplus \fg_2 \oplus \fg_3, 
\quad
\fg_1 =\bR X_1\oplus\bR X_2, \
\fg_2 =\bR X_3,\
\fg_3 =\bR X_4,
$$
where the only non-trivial relations on the canonical basis $\lbrace X_1,X_2,X_3, X_4\rbrace$ are given by
$$
[X_1,X_2]=X_3\ , \ [X_1,X_3]=tX_4\ . 
$$
Then $\fg_t$ is nilpotent, and we denote by   $G_t$ its connected simply connected nilpotent Lie group. 
All the Lie algebra $\fg_t$ and corresponding groups $G_t$, $t>0$, are isomorphic.

\smallskip

We denote by $\theta^1,\theta^2,\theta^3,\theta^4$ 
the basis of $\fg_g^\ast$ dual to the canonical basis $\lbrace X_1,X_2,X_3, X_4\rbrace$.
We have the same decomposition by weights for $\Hom(\bigwedge^\bullet \fg_t,\bR)$:
\begin{align*}
	k=0&\quad 1 \ (w=0)\\
	k=1&\quad \theta^1,\theta^2 \ (w=1)| \ \theta^3 \ (w=2) |\ \theta^4 \ (w=3)\\
	k=2&\quad \theta^{1,2}\ (w=2)| \ \theta^{1,3},\theta^{2,3} \ (w=3)| \ \theta^{1,4},\theta^{2,4} \ (w=4)| \ \theta^{3,4} \ (w=5)\\
	k=3& \quad \theta^{1,2,3}\ (w=4)| \ \theta^{1,2,4}\ (w=5)|\ \theta^{1,3,4},\theta^{2,3,4} \ (w=6)\\
	k=4&\quad \theta^{1,2,3,4} \ (w=7)\,.
\end{align*}
In this diagram, $k$ is the degree of the forms, the numbers in the parentheses indicate the weights of the forms, and the bar signifies a change of weight within each degree.

Throughout this section, we will use the forms $\theta^I$ in the order in the line for $k$ above as the chosen ordered basis of $\Omega^k G$ seen as the $C^\infty (G)$-module generated by $\Hom(\bigwedge^k \fg_t,\bR)$.
For instance, the ordered basis for $\Omega^2 G$ is $\{\theta^{1,2},\theta^{1,3},\theta^{2,3},
\theta^{1,4},\theta^{2,4},  \theta^{3,4} \}$.

Note that with the convention of Section \ref{subsec_weight}, 
$N_0=4$. 

\subsection{The differentials $d$ and 
$d_0$}
\label{subsec_Matd0}

Let us  compute the matrix representations of 
the exterior derivative $d$. 
On $\Omega^0 G$, it is given by
$$
\Mat (d^{(0)})={\tiny \begin{pmatrix}X_1\\ X_2\\ X_3 \\ X_4\end{pmatrix}}.
$$
By the Leibniz rule, we obtain $d^{(k)}$  on $\Omega^k G$, $k>0$, and deduce the following matrix realisations:
\begin{align*}
 \Mat (d^{(1)})&= {\tiny \begin{pmatrix}-X_2 & X_1 & -1 & 0 \\ -X_3 & 0 & X_1 & -t \\ 0 & -X_3 & X_2 & 0 \\ -X_4 & 0 & 0 & X_1\\ 0 & -X_4 & 0 & X_2 \\ 0 & 0 & -X_4 & X_3\end{pmatrix}},\\ 
  \Mat (d^{(2)})&= {\tiny\begin{pmatrix}X_3 & -X_2& X_1 & 0 & -t & 0 \\ X_4& 0 & 0 & -X_2& X_1 & -1\\ 0 & X_4 & 0 & -X_3 & 0 & X_1\\ 0 & 0 & X_4 & 0 & -X_3 & X_2\end{pmatrix}},\\
    \Mat (d^{(3)})&= {\tiny\begin{pmatrix}-X_4 & X_3 & -X_2 & X_1\end{pmatrix}}.
\end{align*}

Let us now describe $d_0$, or equivalently 
$d^{(k)}_0: \Omega^kG\to\Omega^{k+1}G$ for each $k=0,1,2,3$.
We already know that $d_0^{(0)}=0$ and $d_0^{(3)}=0$, whereas the action on the $1$-forms is given by
$$
(d_0^{(1)}\alpha)(V_0,V_1) = -\alpha ([V_0,V_1])\ , 
\quad\alpha\in \Omega^1 G, \ 
V_0,V_1\in \Gamma (TM)\ ,
$$
so that 
\begin{align*}
    d_0^{(1)}\theta^1=d_0^{(1)}\theta^2=0\ ,\quad d_0^{(1)}\theta^3=-\theta^{1,2}\quad  \text{and}\quad d_0^{(1)}\theta^4=-t\theta^{1,3}\,.
\end{align*}
To compute the action of $d_0^{(k)}$ for $k=2$, we apply the Leibniz rule and obtain \begin{align*}
    d_0^{(2)}\theta^{1,2}=d_0^{(2)}\theta^{1,3}=d_0^{(2)}\theta^{2,3}=d_0^{(2)}\theta^{1,4}=0\ ,\quad d_0^{(2)}\theta^{2,4}=-t\theta^{1,2,3}\ ,\quad d_0^{(2)}\theta^{3,4}=-\theta^{1,2,4}\,.
\end{align*}

The action of each $d_0^{(k)}\colon\Omega^kG_t\to\Omega^{k+1}G_t$, 
$k=0,1,2,3,$ can be expressed in block matrix form using the  basis of $\Omega^k G$ described in Section \ref{subsec_GtthetaI}  as follows
\begin{align*}
\Mat (d_0^{(0)}) = 0_{4,1}\ , 
\qquad&
\Mat (d_0^{(3)}) = 0_{1,4}\ ,\\
\Mat (d_0^{(1)}) = 
\begin{pmatrix}
0_{2,2}& \diag(-1,-t) \\
0_{4,2}& 0_{4,2}
\end{pmatrix}\ ,
\qquad&
\Mat (d_0^{(2)}) = 
\begin{pmatrix}
0_{2,4}& \diag(-t,-1) \\
0_{2,4}& 0_{2,2}
\end{pmatrix}\,,
\end{align*}
where $0_{l,k}$ denotes the zero-matrix with $l$ lines and $k$ columns.

\subsection{The operator $\Box_0$}
\label{subsec_Box0Engel}

By definition, we have for $k=0,\ldots,4,$
$$
\Box_0^{(k)}=d_0^{(k-1)}d_0^{(k-1,t)}+d_0^{(k,t)}d_0^{(k)}\colon\Omega^kG_t\to\Omega^kG_t\,, 
$$
which writes matricially as
$$
\Mat (\Box_0^{(k)}) = 
\Mat (d_0^{(k-1)}) 
\big(\Mat (d_0^{(k-1)})\big)^t 
+ \big(\Mat (d_0^{(k)})\big)^t 
\Mat (d_0^{(k)})\,.
$$
Using the computation  for $\Mat (d_0^{(k)})$ in Section \ref{subsec_Matd0}, 
we obtain readily that 
$\Box_0^{(k)}=0$,
and therefore $\Mat (\Box_0^{(k)})=(0)$,
for $k=0,4$, while  for $k=1,2,3$ we have
\begin{align*}
\Mat (\Box_0^{(1)})&=
 \diag (0,0,1,t^2)\,,\\
\Mat (\Box_0^{(2)})&=
\diag (1,t^2,0,0,1,t^2)\,,\\
\Mat (\Box_0^{(3)})&=
    \diag (t^2,1,0,0)\,.
\end{align*}

\subsubsection{The spectral decomposition of $\Box_0$}

By the spectral theorem, 
the operators $\Box_0^{(k)}$
can be expressed as
\begin{align}\label{spectral dec of box0}
    \Box_0^{(k)}=\sum_{j=0}^{m_k}\lambda_j^{(k)}\Pi^{(k)}_{\lambda_j}\ ,\ k=0,\ldots,4\ ,
\end{align}
where $\lambda_0^{(k)},\ldots,\lambda_{m_k}^{(k)}$ are all the distinct eigenvalues, and each $\Pi_{\lambda_i}^{(k)}$ denotes the orthogonal projection onto the eigenspace of $\Box_0^{(k)}$ corresponding to the eigenvalue $\lambda_i^{(k)}$.
We will need to relax this convention with Remark \ref{rem_distinguisht} below.
Indeed, as it is already clear from the matrix realisations of $\Box_0^{(k)}$ above, the dimension of the 1-eigenspace varies when $t=1$ or $t\neq 1$. This means that the spectral resolution of $\Box_0$ depends on the specific realisation of the Engel group $G_t$.

Naturally,  we have for the trivial cases $k=0,4$:
$$
   \Box_0^{(k)}=0\cdot\Pi_0^{(k)}\ , \text{ where }\lambda_0^{(k)}=0\ \text{ and } \Pi_0^{(0)}=\begin{pmatrix}1\end{pmatrix}.
$$
However, for $k=1$ for example, we need to distinguish the two cases $t=1$ and $t\neq 1$.
If $t\neq 1$, then 
$\Box_0^{(1)}=0\cdot\Pi_0^{(1)}+1\cdot\Pi_1^{(1)}+t^2\cdot\Pi_{t}^{(1)}$ has three distinct eigenvalues 
$$
\lambda_0^{(1)}=0, \ \lambda_1^{(1)}=1, 
 \lambda_2^{(1)}=t^2, 
 $$
and associated projectors
 $$
 \Mat(\Pi_{\lambda_0^{(1)}}^{(1)})=\diag (1,1,0,0), 
 $$
 and 
 \begin{equation}
 \label{eq_Pit11}
 \Mat(\Pi_{\lambda_1^{(1)}}^{(1)})=\diag (0,0,1,0), 
 \quad
 \Mat(\Pi_{\lambda_2^{(1)}}^{(1)})=\diag (0,0,0,1).
 \end{equation}
If $t=1$, then 
$\Box_0^{(1)}=0\cdot\Pi_0^{(1)}+1\cdot\Pi_1^{(1)}$ has two distinct eigenvalues $\lambda_0^{(1)}=0$, $\lambda_1^{(1)}=1$
, with associated projectors
 $$
 \Mat(\Pi_0^{(1)})=\diag (1,1,0,0), 
 \qquad
 \Mat(\Pi_1^{(1)})=\diag (0,0,1,1) . 
 $$

\begin{remark}
\label{rem_distinguisht}
In order to avoid having to distinguish between the two cases $t=1$ and $t\neq 1$, we allow for the 
eigenvalues in \eqref{spectral dec of box0} to not be necessarily distinct. Consequently, 
the spectral projection of $\Box_0^{(k)}$ for the eigenvalue $\lambda$ may decompose as 
 $\sum_{j: \lambda_j=\lambda}\Pi_{\lambda_j^{(k)}}^{(k)}$.
With this convention, in the case $t=1$, the spectral projector for $\lambda=1$ is  $\Pi_{\lambda_1^{(1)}}^{(1)} + \Pi_{\lambda_2^{(1)}}^{(1)}$, where $\Pi_{\lambda_1^{(1)}}^{(1)}$ and $ \Pi_{\lambda_2^{(1)}}^{(1)}$ are given matricially via \eqref{eq_Pit11}.
\end{remark}

For $k=2$, we have 
$$
\Box_0^{(2)}=\lambda_0^{(2)} \cdot\Pi_{\lambda_0^{(2)}}^{(2)}
+\lambda_1^{(2)}\cdot\Pi_{\lambda_1^{(2)}}^{(2)}
+\lambda_2^{(2)}\cdot\Pi_{\lambda_2^{(2)}}^{(2)}
\quad\mbox{with}\quad \lambda_0^{(2)}=0, \ \lambda_1^{(2)}=1, \ \lambda_2^{(2)}=t^2,
$$
and 
\begin{align*}
     \Mat(\Pi_{\lambda_0^{(2)}}^{(2)})&=\diag (0,0,1,1,0,0)\,,\\
     \Mat(\Pi_{\lambda_1^{(2)}}^{(2)})&=\diag (1,0,0,0,0,1)\,,\\
     \Mat(\Pi_{\lambda_2^{(2)}}^{(2)})&=\diag (0,1,0,0,1,0)\,.
\end{align*}

For $k=3$, we have
$$
\Box_0^{(3)}=\lambda_0^{(3)} \cdot\Pi_{\lambda_0^{(3)}}^{(3)}
+\lambda_1^{(3)}\cdot\Pi_{\lambda_1^{(3)}}^{(3)}
+\lambda_2^{(3)}\cdot\Pi_{\lambda_2^{(3)}}^{(3)}
\quad\mbox{with}\quad \lambda_0^{(3)}=0, \ \lambda_1^{(3)}=1, \ \lambda_2^{(3)}=t^2,
$$
and 
\begin{align*}
     \Mat(\Pi_{\lambda_0^{(3)}}^{(3)})&=\diag (0,0,1,1)\,,\\
     \Mat(\Pi_{\lambda_1^{(3)}}^{(3)})&=\diag (0,1,0,0)\,,\\
     \Mat(\Pi_{\lambda_2^{(3)}}^{(3)})&=\diag (1,0,0,0)\,.
\end{align*}

As the eigenvalues $\lambda_j^{(k)}$ of $\Box_0^{(k)}$ are the same for $k=1,2,3$, we will allow ourselves to remove their index $k$:
$$ \lambda_0:=\lambda_0^{(k)}=0,\qquad
\lambda_1:=\lambda^{(k)}_1=1, \qquad \lambda_2:=\lambda_2^{(k)}=t^2. 
$$

\subsubsection{The resolvent of $\Box_0$}
With the above conventions, the resolvent  operator $(z-\Box_0)^{-1}$ of $\Box_0$ is given 
for $k=0,4$ by $\big(z-\Box_0^{(k)}\big)^{-1}=\frac{1}{z}\Pi_0^{(k)} = \frac 1z (1)$, 
while for $k=1,2,3$
\begin{equation}
\label{eq_Box0resPi}
\big(z-\Box_0^{(k)}\big)^{-1}=\frac 1z \Pi_0^{(1)}+\frac1{z-1}\Pi_{\lambda_1}^{(k)}+ \frac1{z-t^2} \Pi_{\lambda_2}^{(k)}\,.
\end{equation}

\subsection{The operators $\Box$, $P$ and $L$}

\subsubsection{The operator $\Box$} 

By definition, for each $k=0,\ldots,4$, we have 
$$\Box^{(k)}=d^{(k-1)}d_0^{(k-1,t)}+d_0^{(k,t)}d^{(k)}\colon \Omega^kG_t\to\Omega^{k+1}G_t\,,
$$
so 
$$
\Mat(\Box^{(k)})=\Mat(d^{(k-1)})(\Mat \big(d_0^{(k-1)})\big)^t+\big(\Mat (d_0^{(k)})\big)^t \Mat (d^{(k)})\,.
$$

We have
$\Box^{(k)}=0$
and 
$\Mat(\Box^{(k)})=(0)$ for $k=0,4$.
Using the matrix realisations for $d_0$ and $d$ computed in Section \ref{subsec_Matd0}, 
we compute for $k=1,2,3$:
\begin{align*}
 \Mat(\Box^{(1)})&= {\tiny \begin{pmatrix}0 & 0 & 0 & 0 \\ 0 & 0 & 0 & 0 \\ X_2 & -X_1 & 1 & 0 \\tX_3 & 0 & -tX_1 & t^2 \end{pmatrix}}, \\
  \Mat(\Box^{(2)})&= 
{\tiny \begin{pmatrix}1 &0&0 & 0 & 0& 0 \\
    -X_1 & t^2 & 0 & 0 & 0 & 0 \\ -X_2 & 0 & 0 & 0 & 0 & 0 \\ 0 & -tX_1 & 0 & 0 & 0 & 0 \\ -tX_3 & 0 & -tX_1 & 0& t^2 & 0 \\ 0 & -tX_3 & 0 & X_2 & -X_1 & 1\end{pmatrix}},\\
      \Mat(\Box^{(3)})&= {\tiny \begin{pmatrix}t^2 & 0 & 0 & 0 \\ -tX_1 & 1 & 0 & 0 \\ 0 & -X_1 & 0 & 0 \\ tX_3 & -X_2 & 0 & 0 \end{pmatrix}}.
\end{align*}

\subsubsection{The operator $B=\Box-\Box_0$}
We deduce the matrix representations 
$\Mat(B^{(k)})=\Mat(\Box^{(k)})-\Mat(\Box_0^{k})$
for the operator 
$B^{(k)}=\Box^{(k)}-\Box_0^{k}:\Omega^k G\to \Omega^k G$, $k=0,\ldots,4$.
We have $B^{(k)}=0$, and so $\Mat(B^{(k)})=(0)$, 
for $k=0,4$,
while for $k=1,2,3$
\begin{align*}
 \Mat(B^{(1)})
 &={\tiny \begin{pmatrix}0 & 0 & 0 & 0 \\ 0 & 0 & 0 & 0 \\ X_2 & -X_1 & 0 & 0 \\tX_3 & 0 & -tX_1 & 0 \end{pmatrix}}\,,\\
  \Mat(B^{(2)})
  &={\tiny 
\begin{pmatrix}0 &0&0 & 0 & 0& 0 \\
    -X_1 & 0 & 0 & 0 & 0 & 0 \\ -X_2 & 0 & 0 & 0 & 0 & 0 \\ 0 & -tX_1 & 0 & 0 & 0 & 0 \\ -tX_3 & 0 & -tX_1 & 0& 0 & 0 \\ 0 & -tX_3 & 0 & X_2 & -X_1 & 0 \end{pmatrix}}\,,\\
   \Mat(B^{(3)})
  &= {\tiny 
\begin{pmatrix}0 & 0 & 0 & 0 \\ -tX_1 & 0 & 0 & 0 \\ 0 & -X_1 & 0 & 0 \\ tX_3 & -X_2 & 0 & 0 \end{pmatrix}}\,.
\end{align*}

\subsubsection{The resolvent of $\Box$}

Let us describe  the resolvent  $(z-\Box)^{-1}$ of $\Box$.
Although only the terms containing $\frac 1z$ will occur in our construction, 
we will give  the full computations of $(z-\Box)^{-1}$.

We recall from  Section \ref{subsec_BoxPL} that an application of  Lemma \ref{lem_DHL4.3}
shows that $z-\Box$ is invertible, and that  its inverse is given by the finite sum \eqref{eq_resolventBox} involving $B$ and the resolvent of $\Box_0$.
For $k=0,4$, this gives $(z-\Box^{(k)})^{-1} = z^{-1}$.
For $k=1,2,3$, together with the formula in \eqref{eq_Box0resPi} for the resolvent of $\Box_0$, we obtain
\begin{align*}
  \big(z-\Box^{(1)}\big)^{-1}
  &=\frac{\Pi_0^{(1)}}{z}+\frac{\Pi_{\lambda_1}^{(1)}}{z-1}+\frac{\Pi_{\lambda_2}^{(1)}}{z-t^2}+\frac{\Pi_{\lambda_1}^{(1)}B^{(1)}\Pi_0^{(1)}}{z(z-1)}+\frac{\Pi_{\lambda_2}^{(1)}B^{(1)}\Pi_0^{(1)}}{z(z-t^2)}+\frac{\Pi_{\lambda_2}^{(1)}B^{(1)}\Pi_{\lambda_1}^{(1)}}{(z-1)(z-t^2)}+\\&\qquad+\frac{\Pi_{\lambda_2}^{(1)}B^{(1)}\Pi_{\lambda_1}^{(1)}B^{(1)}\Pi_0^{(1)}}{z(z-1)(z-t^2)}\,,
  \\
  \big(z-\Box^{(2)}\big)^{-1}&=\frac{\Pi_0^{(2)}}{z}+\frac{\Pi_{\lambda_1}^{(2)}}{z-1}+\frac{\Pi_{\lambda_2}^{(2)}}{z-t^2}+\frac{\Pi_0^{(2)}B^{(2)}\Pi_{\lambda_1}^{(2)}+\Pi_{\lambda_1}^{(2)}B^{(2)}\Pi_0^{(2)}}{z(z-1)}+\\&\qquad+\frac{\Pi_0^{(2)}B^{(2)}\Pi_{\lambda_2}^{(2)}+\Pi_{\lambda_2}^{(2)}B^{(2)}\Pi_0^{(2)}}{z(z-t^2)}+\frac{\Pi_{\lambda_1}^{(2)}B^{(2)}\Pi_{\lambda_2}^{(2)}+\Pi_{\lambda_2}^{(2)}B^{(2)}\Pi_{\lambda_1}^{(2)}}{(z-1)(z-t^2)}+\\&\qquad+\frac{\Pi_{\lambda_1}^{(2)}B^{(2)}\Pi_{\lambda_2}^{(2)}B^{(2)}\Pi_{\lambda_1}^{(2)}}{(z-1)^2(z-t^2)}+\frac{\Pi_0^{(2)}B^{(2)}\Pi_{\lambda_2}^{(2)}B^{(2)}\Pi_{\lambda_1}^{(2)}}{z(z-1)(z-t^2)}+
\\&\qquad+\frac{\Pi_{\lambda_1}^{(2)}\big(B^{(2)}\Pi_0^{(2)}B^{(2)}\Pi_{\lambda_2}^{(2)}+B^{(2)}\Pi_{\lambda_2}^{(2)}B^{(2)}\Pi_0^{(2)}\big)+\Pi_{\lambda_2}^{(2)}B^{(2)}\Pi_0^{(2)}B^{(2)}\Pi_{\lambda_1}^{(2)}}{z(z-1)(z-t^2)}+\\&\qquad+\frac{\Pi_{\lambda_1}^{(2)}\big(B^{(2)}\Pi_0^{(2)}B^{(2)}\Pi_{\lambda_2}^{(2)}B^{(2)}\Pi_{\lambda_2}^{(2)}+B^{(2)}\Pi_{\lambda_1}^{(2)}B^{(2)}\Pi_0^{(2)}B^{(2)}\Pi_{\lambda_1}^{(2)}\big)}{z(z-1)^2(z-t^2)}\,,\\
(z-\Box^{(3)})^{-1}
&=\frac{\Pi_0^{(3)}}{z}+\frac{\Pi_{\lambda_1}^{(3)}}{z-1}+\frac{\Pi_{\lambda_2}^{(3)}}{z-t^2}+\frac{\Pi_0^{(3)}B^{(3)}\Pi_{\lambda_1}^{(3)}}{z(z-1)}+\frac{\Pi_0^{(3)}B^{(3)}\Pi_{\lambda_2}^{(3)}}{z(z-t^2)}+\\&\qquad+\frac{\Pi_{\lambda_1}^{(3)}B^{(3)}\Pi_{\lambda_2}^{(3)}}{(z-1)(z-t^2)}+\frac{\Pi_0^{(3)}B^{(3)}\Pi_{\lambda_1}^{(3)}B^{(3)}\Pi_{\lambda_2}^{(3)}}{z(z-1)(z-t^2)}.
    \end{align*}

\subsubsection{The projection $P$}

By definition, $P=\frac{1}{2\pi i}\oint_{\vert z\vert =\varepsilon}\big(z-\Box\big)^{-1}dz$ is the Cauchy integral about 0 of the resolvent of $\Box$.
The above computations for the resolvent of $\Box_0$ yield $P^{(k)} = (1)$ for $k=0,4$ and 
 for $k=1,2,3$: 
\begin{align*}
    &P^{(1)}=
    \Pi_0^{(1)}-\Pi_{\lambda_1}^{(1)}B^{(1)}\Pi_0^{(1)}+\frac{\Pi_{\lambda_2}^{(1)}B^{(1)}\Pi_{\lambda_1}^{(1)}B^{(1)}\Pi_0^{(1)}-\Pi_{\lambda_2}^{(1)}B^{(1)}\Pi_0^{(1)}}{t^2}\,\\
    &P^{(2)}=\Pi_0^{(2)}-\Pi_0^{(2)}B^{(2)}\Pi_{\lambda_1}^{(2)}-\Pi_{\lambda_1}^{(2)}B^{(2)}\Pi_0^{(2)}-\frac{\Pi_0^{(2)}B^{(2)}\Pi_{\lambda_2}^{(2)}+\Pi_{\lambda_2}^{(2)}B^{(2)}\Pi_0^{(2)}}{t^2}+\\&\quad+\frac{\Pi_0^{(2)}B^{(2)}\Pi_{\lambda_2}^{(2)}B^{(2)}\Pi_{\lambda_1}^{(2)}+\Pi_{\lambda_1}^{(2)}\big(B^{(2)}\Pi_0^{(2)}B^{(2)}\Pi_{\lambda_2}+B^{(2)}\Pi_{\lambda_2}^{(2)}B^{(2)}\Pi_0^{(2)}\big)}{t^2}+\\&\quad+\frac{\Pi_{\lambda_2}^{(2)}B^{(2)}\Pi_0^{(2)}B^{(2)}\Pi_{\lambda_1}^{(2)}-\Pi_{\lambda_1}^{(2)}\big(B^{(2)}\Pi_0^{(2)}B^{(2)}\Pi_{\lambda_2}^{(2)}B^{(2)}\Pi_{\lambda_2}^{(2)}+B^{(2)}\Pi_{\lambda_1}^{(2)}B^{(2)}\Pi_0^{(2)}B^{(2)}\Pi_{\lambda_1}^{(2)}\big)}{t^2}\,,\\
    &P^{(3)}=\Pi_0^{(3)}-\Pi_0^{(3)}B^{(3)}\Pi_{\lambda_1}^{(3)}-\frac{\Pi_0^{(3)}B^{(3)}\Pi_{\lambda_2}^{(3)}}{t^2}+\frac{\Pi_0^{(3)}B^{(3)}\Pi_{\lambda_1}^{(3)}B^{(3)}\Pi_{\lambda_2}^{(3)}}{t^2}\,.
\end{align*}
Using the matrix representations of the $\Pi_{\lambda_i}^{(k)}$ from Section \ref{subsec_Box0Engel},  we obtain:
\begin{align*}
\Mat (P^{(1)})
&={\tiny \begin{pmatrix}1 & 0 & 0 & 0\\0 &1 & 0 & 0 \\-X_2 & X_1 & 0 & 0\\-(X_1X_2+X_3)/t & X_1^2/t & 0 & 0 \end{pmatrix}}\,,
\\    
\Mat (P^{(2)})
&={\tiny \begin{pmatrix}0 &0&0 & 0 & 0& 0 \\
    0 & 0 & 0 & 0 & 0 & 0 \\ 
    X_2 & 0 & 1 & 0 & 0 & 0 \\ 
    X_1^2/t & X_1/t & 0 & 1 & 0 & 0 \\ 
    X_1X_2/t & 0 & X_1/t & 0& 0 & 0 \\ (X_1^2X_2-X_2X_1^2)/t & -X_2X_1/t & X_1^2/t & -X_2 & 0 & 0\end{pmatrix}},\\    
\Mat (P^{(3)})
&={\tiny \begin{pmatrix}0 & 0 & 0 & 0\\ 0 & 0 & 0 & 0\\ X_1^2/t & X_1 & 1 & 0 \\(X_2X_1- X_3)/t & X_2 & 0 & 1\end{pmatrix}}.
\end{align*}

\subsubsection{The operator $L$}
The differential operator $L^{(k)}\colon\Omega^kG\to\Omega^{k}G$ is given by
\begin{align*}
    L^{(k)}=P^{(k)}\Pi_0^{(k)}+(\id-P^{(k)})(\id-\Pi_0^{(k)})\,, 
\end{align*}
so 
\begin{align*}
    \Mat(L^{(k)})=\Mat(P^{(k)})\Mat(\Pi_0^{(k)})+(\id-\Mat(P^{(k)}))(\id-\Mat(\Pi_0^{(k)}))\,.
\end{align*}
We have $\Mat (L^{(k)}) = (1)$ for $k=0,4$, 
while for $k=1,2,3$
\begin{align*}
  L^{(1)}&=
  {\tiny \begin{pmatrix}1 & 0 & 0 & 0 \\ 0 & 1 & 0 & 0\\ -X_2 & X_1 & 1 & 0 \\ -(X_1X_2+X_3)/t & X_1^2/t & 0& 1  \end{pmatrix}}\,,  \\
  L^{(2)}&=
  {\tiny \begin{pmatrix}1 &0&0 & 0 & 0& 0 \\
    0 & 1 & 0 & 0 & 0 & 0 \\ -X_2 & 0 & 1 & 0 & 0 & 0 \\ -X_1^2/t & -X_1/t & 0 & 1 & 0 & 0 \\ -X_1X_2/t & 0 & X_1/t & 0& 1 & 0 \\ (X_2X_1^2-X_1^2X_2)/t & X_2X_1/t & X_1^2/t & -X_2 & 0 & 1\end{pmatrix}}\,,
  \\L^{(3)}&={\tiny \begin{pmatrix}1 & 0 & 0 & 0 \\ 0 & 1 & 0 & 0\\ -X_1^2/t & -X_1 & 1 & 0 \\ (X_3-X_2X_1)/t & -X_2 & 0 &1  \end{pmatrix}}\,.
\end{align*}

\subsubsection{The operator $L^{-1}$}

The inverse  operator $L^{-1}$ is given by 
$$
\big(L^{(k)}\big)^{-1}=\sum_{j=0}^{N_0-1}\big(\id-L^{(k)}\big)^j.
$$
and matricially, 
$$
\big(\Mat(L^{(k)})\big)^{-1}=\Mat\left(\big(L^{(k)}\big)^{-1}\right)=\sum_{j=0}^{N_0-1}\big(\id- \Mat(L^{(k)})\big)^j\,.
$$
Hence, we have $\big(\Mat(L^{(k)})\big)^{-1}=(1)$ for $k=0,4$ while for $k=1,2,3,$ we obtain using the above matrix realisations for $L^{(k)}$:
\begin{align*}
 \big(\Mat(L^{(1)})\big)^{-1}   
 &={\tiny \begin{pmatrix}1 & 0 & 0 & 0 \\ 0 & 1 & 0 & 0\\ X_2 & -X_1 & 1 & 0 \\ (X_1X_2+X_3)/t & -X_1^2/t & 0& 1  \end{pmatrix}}\,, \\
 \big(\Mat(L^{(2)})\big)^{-1}
 &={\tiny \begin{pmatrix}1 &0&0 & 0 & 0& 0 \\
    0 & 1 & 0 & 0 & 0 & 0 \\ X_2 & 0 & 1 & 0 & 0 & 0 \\ X_1^2/t & X_1/t & 0 & 1 & 0 & 0 \\ 0 & 0 & -X_1/t & 0& 1 & 0 \\ 0 & 0 & -X_1^2/t & X_2 & 0 & 1\end{pmatrix}}\,,\\
 \big(\Mat(L^{(3)})\big)^{-1}
 &={\tiny \begin{pmatrix}1 & 0 & 0 & 0 \\ 0 & 1 & 0 & 0\\ X_1^2/t & X_1 & 1 & 0 \\ (X_2X_1-X_3)/t & X_2 & 0& 1  \end{pmatrix}}\,.
\end{align*}

\subsection{The operators $D$ and $C$}

\subsubsection{The operator $L^{-1}dL$}
In order to determine $D$ and $C$, we first compute $L^{-1}dL$, or more precisely
$$
(L^{(k+1)})^{-1}d^{(k)}L^{(k)}\colon\Omega^kG\to\Omega^{k+1}G\,, \quad k=0,1,2,3,4\,.
$$
Its matrix realisation is
$$
\Mat\big((L^{(k+1)})^{-1}d^{(k)}L^{(k)}\big)
\big(\Mat (L^{(k+1)})\big)^{-1} \cdot\Mat(d^{(k)})\cdot\Mat(L^{(k)}),
$$
 and from the previous computations for the matrix representations of $L,L^{-1},d,$ we obtain:
\begin{align*}
 \Mat\big((L^{(1)})^{-1}d^{(0)}L^{(0)}\big)&= {\tiny \begin{pmatrix}X_1\\ X_2\\ 0 \\ 0\end{pmatrix}}\,,\\
\Mat\big( (L^{(2)})^{-1}d^{(1)}L^{(1)}\big)&=
{\tiny  \begin{pmatrix}0& 0 & -1 & 0 \\ 0 & 0 & X_1 & -t\\ -X_2^2& X_2X_1-X_3& 0& 0\\ -X_4-(X_1^2X_2+X_1X_3)/t& X_1^3/t& 0 & 0\\0&0&-{X_1X_2}/t & X_2\\ 0& 0 & -X_4-X_1^2X_2/t& X_3+X_2X_1\end{pmatrix}},
 \\
\Mat\big( (L^{(3)})^{-1}d^{(2)}L^{(2)}\big)&=
{\tiny \begin{pmatrix}X_3& -X_2 & 0 & 0 & -t & 0\\ 
    X_4& 0 & 0 & 0 & X_1 & -1\\ 0 & 0 & {X_1^3}/t& -X_3-X_1X_2& 0 &0\\0 & 0 & ({X_1^2X_2-3X_3X_1})/t& -X_2^2 & 0 & 0 \end{pmatrix}}\,,
    \\
\Mat\big((L^{(4)})^{-1}d^{(3)}L^{(3)}\big)&=
{\tiny \begin{pmatrix}
    0 & 0& -X_2&X_1
    \end{pmatrix}}\,.
\end{align*}

\subsubsection{The operator $D$}
We can now compute explicitly the operator $D=\Pi_0L^{-1}dL\Pi_0$ in matrix terms:
\begin{align*}
\Mat (D^{(0)})&={\tiny \begin{pmatrix}
    X_1\\ X_2 \\ 0 \\ 0 
    \end{pmatrix}}\,,\\
 \Mat (D^{(1)})&={\tiny \begin{pmatrix}0& 0 & 0 & 0 \\ 0 & 0 & 0 & 0\\ -X_2^2& -X_3+X_2X_1& 0& 0\\ -X_4-({X_1^2X_2+X_1X_3})/t& {X_1^3}/t& 0 & 0\\0&0&0 & 0\\ 0& 0 & 0& 0\end{pmatrix}}\,, \\
 \Mat (D^{(2)})&={\tiny \begin{pmatrix}0 & 0 & 0 & 0 & 0 & 0 \\ 0 & 0 & 0 & 0 & 0 & 0 \\0 & 0 & {X_1^3}/t & -X_3-X_1X_2 & 0 & 0 \\0 & 0 & ({X_1^2X_2-3X_3X_1})/t & -X_2^2 & 0 & 0 \end{pmatrix}}\,, \\
 \Mat (D^{(3)})&={\tiny \begin{pmatrix}0 & 0 & -X_2 & X_1 \end{pmatrix}}\,.
\end{align*}

\subsubsection{The operator $C$}
Similarly, we obtain for $C =(\id- \Pi_0)L^{-1}dL(\id-\Pi_0)$:
\begin{align*}
\Mat (C^{(0)})&= 0_{4,1}\,,\\
\Mat (C^{(1)})&={\tiny \begin{pmatrix}0& 0 & -1 & 0 \\ 0 & 0 & X_1 & -t\\ 0& 0& 0& 0\\ 0& 0& 0 & 0\\0&0&-{X_1X_2}/t & X_2\\ 0& 0 & -X_4-{X_1^2X_2}/t& X_1X_2\end{pmatrix}}\,,\\
\Mat (C^{(2)})&={\tiny \begin{pmatrix}X_3 & -X_2 & 0 & 0 & -t & 0 \\ X_4 & 0 & 0 & 0 & X_1 & -1 \\0 & 0 & 0 & 0 & 0 & 0 \\0 & 0 & 0 & 0 & 0 & 0 \end{pmatrix}}\,,\\
\Mat (C^{(3)})&= 0_{1,4}\,.
\end{align*}

	\bibliographystyle{alpha.bst}

\bibliography{bibli}


\end{document}